\documentclass[11pt,leqno]
{amsart}
\usepackage{amsmath,epsfig,graphicx,color}
\definecolor{darkblue}{rgb}{.2, 0.2,.8}
\definecolor{darkgreen}{rgb}{0,0.5,0.3}
\definecolor{darkred}{rgb}{.8, .1,.1}

\newcommand{\levy}{L\'evy}

\newcommand{\ex}{{\rm e}\,}

\newcommand{\asy}{asymptotic}

\newcommand {\spp} {^\prime}
\newcommand{\ts}{time series}
\newcommand{\tsa}{\ts\ analysis}

\newtheorem{lemma}{Lemma}[section]

\newtheorem{theorem}[lemma]{Theorem}

\newtheorem{proposition}[lemma]{Proposition}
\newtheorem{definition}[lemma]{Definition}
\newtheorem{corollary}[lemma]{Corollary}
\newtheorem{example}[lemma]{Example}
\newtheorem{exercise}[lemma]{Exercise}
\newtheorem{remark}[lemma]{Remark}
\newtheorem{fig}[lemma]{Figure}
\newtheorem{tab}[lemma]{Table}

\newcommand{\RV}{{\bf RV}}

\newcommand{\bth}{\begin{theorem}}
\newcommand{\ethe}{\end{theorem}}
\newcommand{\sv}{stochastic volatility}
\newcommand{\bre}{\begin{remark}\em }
\newcommand{\ere}{\end{remark}}

\newcommand{\ble}{\begin{lemma}}
\newcommand{\ele}{\end{lemma}}
\newcommand{\sre}{stochastic recurrence equation}

\newcommand{\bde}{\begin{definition}}
\newcommand{\ede}{\end{definition}}
\newcommand{\bco}{\begin{corollary}}
\newcommand{\eco}{\end{corollary}}

\newcommand{\bpr}{\begin{proposition}}
\newcommand{\epr}{\end{proposition}}

\newcommand{\bexer}{\begin{exercise}}
\newcommand{\eexer}{\end{exercise}}

\newcommand{\bexam}{\begin{example}}
\newcommand{\eexam}{\end{example}}

\newcommand{\bfi}{\begin{fig}}
\newcommand{\efi}{\end{fig}}

\newcommand{\btab}{\begin{tab}}
\newcommand{\etab}{\end{tab}}

\newcommand{\rv}{random variable}

\newcommand{\MDA}{{\rm MDA}}
\newcommand{\sign}{{\rm sign}}

\newcommand{\var}{{\rm var}}

\newcommand{\bfTh}{\mbox{\boldmath$\Theta$}}

\newcommand{\rhs}{right-hand side}
\newcommand{\df}{distribution function}

\newcommand{\beao}{\begin{eqnarray*}}
\newcommand{\eeao}{\end{eqnarray*}\noindent}

\newcommand{\beam}{\begin{eqnarray}}
\newcommand{\eeam}{\end{eqnarray}\noindent}

\newcommand{\beqq}{\begin{equation}}
\newcommand{\eeqq}{\end{equation}\noindent}

\newcommand{\bce}{\begin{center}}
\newcommand{\ece}{\end{center}}

\newcommand{\barr}{\begin{array}}
\newcommand{\earr}{\end{array}}

\newcommand{\stp}{\stackrel{\P}{\rightarrow}}
\newcommand{\std}{\stackrel{d}{\rightarrow}}

\newcommand{\stw}{\stackrel{w}{\rightarrow}}

\newcommand{\eqd}{\stackrel{d}{=}}

\newcommand{\vague}{\stackrel{\lower0.2ex\hbox{$\scriptscriptstyle
                    \it{v} $}}{\rightarrow}}
\newcommand{\weak}{\stackrel{\lower0.2ex\hbox{$\scriptscriptstyle
                    \it{w} $}}{\rightarrow}}
\newcommand{\what}{\stackrel{\lower0.2ex\hbox{$\scriptscriptstyle
                    \it{\hat{w}} $}}{\rightarrow}}

\newcommand{\bdis}{\begin{displaymath}}
\newcommand{\edis}{\end{displaymath}\noindent}

\newcommand{\R}{\mathbb{R}}

\newcommand{\nto}{n\to\infty}

\newcommand{\xto}{x\to\infty}

\newcommand{\ov}{\overline}
\newcommand{\wt}{\widetilde}
\newcommand{\wh}{\widehat}
\newcommand{\vep}{\varepsilon}

\newcommand{\regvary}{regularly varying}
\newcommand{\slvary}{slowly varying}
\newcommand{\regvar}{regular variation}

\newcommand{\bbr}{{\mathbb R}}

\newcommand{\bbz}{{\mathbb Z}}
\newcommand{\bbn}{{\mathbb N}}

\newcommand{\con}{convergence}

\newcommand{\st}{such that}
\newcommand{\fif}{if and only if}

\newcommand{\fct}{function}

\newcommand{\ds}{distribution}

\newcommand{\rep}{representation}

\newcommand{\seq}{sequence}

\newcommand{\ins}{insurance}

\newcommand{\pro}{probabilit}

\newcommand{\ms}{measure}
\newcommand{\mgf}{moment generating function}
\newcommand{\ld}{large deviation}

\newcommand{\bfX}{{\bf X}}

\newcommand{\bfY}{{\bf Y}}

\newcommand{\E }{{\mathbb E}}
\renewcommand{\P }{{\mathbb P}}

\newcommand{\1}{{\mathbf 1}}

\allowdisplaybreaks

\begin{document}
\today
\bibliographystyle{plain}
\title[Precise large deviations for dependent subexponential variables]{Precise large deviations for
dependent subexponential variables}
\author[T. Mikosch]{Thomas Mikosch}
\address{Department  of Mathematics,
University of Copenhagen,
Universitetsparken 5,
DK-2100 Copenhagen,
Denmark}
\email{mikosch@math.ku.dk}
\author[I. Rodionov]{Igor Rodionov}
\address{V.A. Trapeznikov Institute of Control Sciences of the Russian Academy of Sciences, Profsoyuznaya ulitsa 65, 117997, Moscow, Russia}
\email{vecsell@gmail.com}

\begin{abstract}
In this paper we study precise large deviations for
the partial sums of a stationary \seq\ with a subexponential marginal
\ds . Our main focus is on  \ds s which  either have a
\regvary\ or a lognormal-type tail. We apply the results to prove limit
theory for the maxima of the entries large sample covariance matrices.
\end{abstract}
\maketitle

\section{Introduction}\setcounter{equation}{0}
We consider a (strictly) stationary real-valued
\seq\ $(X_t)$ with generic element $X$ and \df\ $F$ with finite first moment.
The corresponding centered partial sums are given by
\beam\label{eq:feb20g}
S_0=0\,,\qquad S_n=X_1+\cdots +X_n-n\,\E[X]\,,\qquad n\ge 1\,.
\eeam
To ease notation we will always assume  that $X$ is centered. We also assume that $F$ is {\em subexponential}.
\subsection{Subexponential \ds s}\label{subsec:subexpdist}
For the moment assume $(X_i)$ are iid.
Following the classical definition of
\cite{cistyakov:1964} (cf. \cite{embrechts:klueppelberg:mikosch:1997}, p.~39),
$F$ is subexponential if $X$ is non-negative
and has the tail-equivalence
property for convolutions, i.e.,
\beam\label{eq:feb20a}
\P(S_n>x)\sim n\,(1- F(x))=n\,\ov F(x)\,,\qquad n\ge 2\,,\qquad \xto \,;
\eeam
we write $F\in\mathcal S_+$.
Here $f(x)\sim g(x)$ for positive \fct s $f,g$ means that $f(x)/g(x)\to1$ as
$\xto$.
In this paper we will consider two-sided
subexponential \df s, i.e., $X^+=X\vee 0$ has a subexponential \ds\ and a {\em tail balance condition} holds
\beam \label{tail-balance}
\lim_{x\to\infty} \dfrac{\P(X> x)}{\P(|X|>x)} = p_{+}, \quad \lim_{x\to\infty} \dfrac{\P(-X> x)}{\P(|X|>x)} = p_{-}
\eeam
for some $p_+>0, p_-\ge 0,$ and we write $\mathcal S$ for this
enlarged class of \ds s. The property \eqref{eq:feb20a} has the interpretation
that $S_n$ and $M_n=\max(X_1,\ldots,X_n)$ are tail-equivalent for every $n$. Therefore it is considered a very natural class  of {\em heavy-tailed \ds s}
which has multiple applications in \ins\ mathematics, telecommunications,
queuing and branching theory. Textbook treatments can be found in \cite{embrechts:klueppelberg:mikosch:1997},
\cite{rolski:schmidli:schmidt:teugels:1999}, and
\cite{asmussen:2003}.
\par
The class $\mathcal S_+$ covers a wide range of tail behaviors from power laws
with certain moments infinite to semi-exponential tails \st\ $X$ has all
moments finite but no \mgf . We will mainly be interested in two sub-classes of
\df s $F\in\mathcal S$:
\begin{itemize}
\item {\bf RV}$(\alpha)$. We say that $X$ and its \ds\ $F$ are \regvary\ with
index $\alpha>0$ $(F\in {\bf RV}(\alpha))$
if $\ov F(x)=1-F(x)=L(x) x^{-\alpha}$ for some \slvary\
\fct\ $L$.
\item {\bf LN}. This class consists of subexponential \ds s $F$ \st\
$\ov F(x)=\exp(-S(x))$ where $S$ is a \slvary\ \fct\ \st\
\beao
S(x)/\log x\to\infty\,,\qquad \xto\,.
\eeao
\end{itemize}
Well-known representatives $F\in{\bf RV}(\alpha)$
with positive tail index $\alpha$ are the Pareto, Burr, student \ds s.
A representative of {\bf LN} is the (standard) lognormal \ds\ with tail
\beam\label{eq:feb20d}
\P(X>x)\sim \dfrac{\ex^{-(\log x)^2/2}}{\sqrt{2\pi} \log x}
=\ex^{-(\log x)^2/2+\log(\sqrt{2\pi} \log x)}\,.
\eeam
An interesting third subclass of $\mathcal S_+$ are the Weibull-type \ds s
with tail $\ov F(x)=\exp(-x^{\alpha}L(x))$
for a \slvary\ \fct\ $L$ and $\alpha\in (0,1)$. Unfortunately, the techniques developed in
this paper fail for these \ds s, see Remark \ref{reWeibull} below.

\subsection{Precise \ld s of subexponential type in the iid case}\label{subsec:ldsubexp}
Early on, it was discovered that the defining property of a subexponential
\ds\ \eqref{eq:feb20a} extends to situations when $n\to\infty$ and
$x=x_n\to\infty$. To be more precise, a relation of the type
\beam\label{eq:feb20c}
\sup_{x>t_n} \Big|\dfrac{\P(S_n>x)}{n\,\ov F(x)}-1\Big|\to 0\,,\qquad \nto\,,
\eeam
holds for a suitable \seq\ $(t_n)$; we call it a {\em separating \seq },
and \eqref{eq:feb20c} a {\em (precise) \ld\ of subexponential type.}
As a matter of fact, \cite{cline:hsing:1998}
discovered that $F\in\mathcal S$ is an ``almost'' necessary and sufficient
condition for \eqref{eq:feb20c} to hold. Pioneering work on \ld s of
type \eqref{eq:feb20c} is due to
\cite{nagaev:1969a,nagaev:1969b,nagaev:1977},  \cite{nagaev:1965,nagaev:1979},  \cite{rozovski:1993}; see also
 \cite{cline:hsing:1998},   \cite{denisov:dieker:shneer:2008}.  Large deviations
for the sample
paths of a \levy\ process and random walks with \regvary\ increments
were considered by
\cite{hult:lindskog:mikosch:samorodnitsky:2005}, \cite{rhee:blanchet:zwart:2019}.
\par
The perhaps best known result in this context is due to
 \cite{nagaev:1979}. For $F\in{\bf RV}(\alpha)$ and $\alpha>2$, assuming $\E[X]=0$ and $\var(X)=1$,
he proved that \eqref{eq:feb20c} holds for $x>t_n=\sqrt{(\alpha-2) n\,\log n}$, while for $x<t_n$ one has
\beam \label{normal-case}
\sup_{x<t_n} \Big|\dfrac{\P(S_n>x)}{\ov \Phi(x/\sqrt{n})}-1\Big|\to 0\,,\qquad \nto\,,
\eeam
where $\Phi$ is the standard normal \df .
\par
Results of the types of \eqref{eq:feb20c} and \eqref{normal-case}
are also valid for various other \ds s in $\mathcal S$. In particular,
the lognormal \ds\ with tail \eqref{eq:feb20d} satisfies \eqref{eq:feb20c}
for $x\gg t_n$ and \eqref{normal-case}
for $x\ll t_n$ where  $t_n=\sqrt{n}\,\log n$ and $x\gg t_n$ means that
$x\ge t_n\,h_n$ for any \seq\ $h_n\to\infty$, and
$x\ll t_n$ is defined correspondingly.  \cite{rozovski:1993}
found that the separating \seq s $(t_n)$ in \eqref{eq:feb20c} and \eqref{normal-case} have to be distinct if $\ov F$ is lighter than the tail of a lognormal
\ds .
\par
Extensions of \ld s of subexponential type to stationary \seq s
only exist in a few cases. \cite{mikosch:samorodnitsky:2000} proved \ld s of subexponential
type for
{\em  \regvary } linear processes driven by iid \regvary\ noise.
The main difference to the iid case is that
the limit of $\P(S_n>x)/(n\ov F(x))$ converges uniformly for $x\gg t_n$
to a constant depending on the coefficients of the linear process and the
tail index of the noise. This fact shows that {\em extremal clustering}
in the $X$-\seq\ causes that exceedances of $S_n$ above high thresholds $x$
appear in clumps and not separated from each other, and the limiting constant
is a \ms\ of the size of these clumps.
Solutions to affine \sre s $X_t=A_tX_{t-1}+B_t$, $t\in\bbz$, for an iid \seq\ $(A_t,B_t)$, $t\in\bbz$, may
have power-law tails $\P(\pm X>x)\sim c_{\pm}x^{-\alpha}$ for some $\alpha>0$
either due to \regvar\ of $B_1$ with index $\alpha$ and $\E[|A_1|^\alpha]<1$
(the so-called Grincevi\v cius-Grey case) or due to the condition
$\E[|A_1|^\alpha]=1$ (the so-called Kesten-Goldie case); see Section 3.4.2 in \cite{buraczewski:damek:mikosch:2016} for an overview.
 \cite{buraczewski:damek:mikosch:zienkiewicz:2013}
proved \ld\ results of subexponential type in the Kesten-Goldie case,
and  \cite{konstantinides:mikosch:2005}
in the Grincevi\v cius-Grey case.
\cite{mikosch:wintenberger:2013,mikosch:wintenberger:2016}
derived \ld\ results for \regvary\ Markov chains and $m$-dependent processes
and applied these results to get bounds for ruin \pro ies.

\subsection{Goals of this paper}
In this paper we aim at proving analogs of the subexponential \ld\
results for a stationary dependent \seq\ $(X_t)$.
In most cases, we have to restrict
ourselves to an $m$-dependent \seq , i.e., the dependence ranges only over
$m$ lags. We work under the heavy-tail assumption $F\in \mathcal S$ which is a
natural condition, as we explained in Section~\ref{subsec:ldsubexp}.
We also have to impose an {\em \asy\ tail independence condition} on the
\ds s of the pairs $(X_0,X_h)$ for $1\le h\le m$. Under the aforementioned
conditions and for $F\in{\bf RV(\alpha)}$ and $F\in{\bf LN}$ we prove results of the type
\eqref{eq:feb20c}. The strong \asy\ tail independence conditions ensure that
\eqref{eq:feb20c} is valid for suitable \seq s $(t_n)$. Based on the
$m$-dependence of $(X_t)$ we make heavy use of the known \ld\
results in the iid case. This is the topic of Section~\ref{sec:main}.
\par
In Section~\ref{subsec:linear} we study subexponential \ld s for
a linear process driven by an iid noise \seq\ with a common subexponential
\ds\ $F$ in the class {\bf LN}. In this case, a result of type \eqref{eq:feb20c}
does in general not hold but the denominator $n\,\ov F(x)$ has to replaced
by $n\,\ov F(x/|m_0|)$ for some number $m_0$ which depends on the
coefficients of the linear process. The proof makes heavy use of the
linear structure and exploits the known \ld\ results for an iid \seq .
We also mention that the \ds\ of $X$ is tail-equivalent to the subexponential
noise \ds .
\par
In Section~\ref{sec:matrix} we show how \ld s of subexponential type
can be applied to determine the limits  of the
maxima of the diagonal or off-diagonal entries of a large sample
covariance matrix with
row-wise dependent entries.

\section{Preliminaries}\setcounter{equation}{0}
\subsection{Maximum domains of attraction}\label{subsec:mda}
Assume that $(X_i)$ is iid with common \ds\ $F$.
\par
The condition $F\in {\bf RV}(\alpha)$ for $\alpha>0$ is equivalent to
membership of $F$ in  the {\em maximum domain of attraction} of the Fr\'echet \ds\ $\Phi_\alpha$
$(F\in\MDA(\Phi_\alpha))$. This means that there exist constants
$a_n>0$ \st\
\beao
\P\big(a_n^{-1} M_n\le x\big)\to \Phi_\alpha(x)=\ex^{-x^{-\alpha}}\,,\qquad x\ge 0\,,\qquad \nto\,.
\eeao
For $F\in{\bf LN}\cap\mathcal S$ we also require that it is a member
of the maximum domain of attraction of the Gumbel \ds\ $\Lambda$
$(F\in\MDA(\Lambda))$, i.e., there exist constants $c_n>0,d_n\in\bbr$
\st\
\beao
\P\big(c_n^{-1}(M_n-d_n)\le x\big)\to \Lambda(x)=
\ex^{-\ex^{-x}}\,,\qquad x\in\bbr\,,\qquad  \nto\,.
\eeao
\par
According to the Pickands-Balkema-de Haan
Theorem (\cite{pickands:1975,balkema:dehaan:1974}, cf. Theorem~3.4.5
in \cite{embrechts:klueppelberg:mikosch:1997}) a \ds\ with
infinite right endpoint $F\in\MDA(\Lambda)$  \fif\
there exists a positive \fct\ $a$ with Lebesgue density $a'$ \st\  $a'(x)\to 0$ as $\xto$ and
\beam\label{gumbel}
\dfrac{\ov F(x+ y\,a(x))}{\ov F(x)} \to \ex^{-y}\,,\qquad \xto\,,\qquad y\in\bbr\,.
\eeam
The auxiliary  \fct\ $a$ can be chosen as the {\em mean-excess \fct\ of $F$}
\beao
a(x) = \int_x^\infty \dfrac{\ov F(y)}{\ov F(x)}\,dy\,,\qquad x>0\,;
\eeao
cf. \cite{resnick:1987}, Proposition~1.9. We have $F\in\MDA(\Lambda)$ \fif\
\beam\label{von_mises} \ov F(x) = c(x)\exp\Big(-\int_z^x \frac{1}{a(t)}dt\Big), \qquad x > z,
\eeam
for some $z$ and $c(x)\to c>0$ as $\xto.$
\subsection{Long-tailed distributions}

A distribution function $F$ is said to be long-tailed
 if
\beao
\lim_{x\to\infty}\dfrac{\ov F(x+y)}{\ov F(x)} = 1, \quad \mbox{for any } y>0.
\eeao
For the properties of long-tailed distributions we refer to \cite{foss:korshunov:zachary:2011}. In particular, $F\in\mathcal{S}$ implies long-tailedness
of $F$; see  Lemma 3.4 in \cite{foss:korshunov:zachary:2011}. Moreover, for each long-tailed distribution $F$ there exists a non-decreasing function $h$ with $h(x) \uparrow \infty$ as $\xto$ such that
\beao
\lim_{x\to\infty}\dfrac{\ov F(x + h(x))}{\ov F(x)} = 1,
\eeao
and $F$ is called {\it $h$-insensitive}. In particular,
$F\in\MDA(\Lambda)$ satisfies \eqref{gumbel} for some auxiliary
\fct\ $a(x)\to\infty$. Hence we can choose $h(x)=o(a(x))$, and if $F\in \MDA(\Phi_{\alpha})$
we can take any \fct\ $h$ with $h(x)=o(x)$ as $\xto$.

\subsection{Condition {\bf (C)}}
We consider a stationary \seq\ $(X_i)$ with mean zero and partial sum process $(S_n)$ given in
\eqref{eq:feb20g}.
In this section we assume that $F\in\MDA(\Lambda)\cap {\bf LN}$. Hence, in particular,
$F\in\mathcal S$, $F$ has infinite right endpoint and
$S(x)=-\log \ov F(x)$ is slowly varying \st\ $S(x)/\log x\to \infty$ as $\xto$. Characterizations
of $\MDA(\Lambda)$ are given in Section~\ref{subsec:mda}.
\par
In what follows, we introduce and discuss a set of conditions
which will be assumed in our main result, Theorem~\ref{thm:main}.
A crucial object in this context is a positive \fct\
$g$ which describes the region $(t_n,\infty)$ where the \ld\ results hold.
\subsection*{Condition {\bf (C)}}
\begin{enumerate}
\item[$\bf C_1$] $g(x)\uparrow \infty$ as $\xto$ and
 there is $C>0$ \st\ for large $x$,
\beam \label{eq:4} g(x)\le C\,x\,/S(x)\,.
\eeam
\item[$\bf C_2$]
There is a \seq\ $t_n\to\infty$ \st\ for any $\delta>0$,
\beam\label{eq:3}
\sup_{x>t_n\delta}\Big|\dfrac{S(x)}{S(g(x))}-1\Big|\to 0\quad\mbox{ and }\quad \dfrac{g(t_n)}{\sqrt{n}} \to \infty\,, \qquad \nto,
\eeam
and for an
iid \seq\ $(X_i')$ with common \ds\ $F$ and partial sums
$S_n'=X_1'+\cdots +X_n'$ we have the \ld\ result
\beam\label{eq:feb20h}
\lim_{\nto}\sup_{x>t_n\delta} \left|\dfrac{\P(S_n'>x)}{n \ov F(x)} - 1\right| = 0\,,\qquad \mbox{ for any $\delta>0$.}
\eeam
\item[$\bf C_3$] $(X_i)$ is $m$-dependent for some $m\ge 1$, and for any $\vep>0$,
\beam\label{eq:1}
\lim_{\xto}\dfrac{\P(|X_0|>\varepsilon g(x)\,,|X_h|>\varepsilon x)}{\ov F(x)}=0\,,\qquad
h=1,\ldots, m\,.
\eeam
\end{enumerate}
\noindent {\bf The size of $g(x)$.} It follows from the monotone density theorem (cf. Theorem 1.7. in
 \cite{bingham:goldie:teugels:1987})
and \eqref{von_mises} that
\beao
\dfrac{a(x)S(x)}{x} \to \infty, \qquad x\to\infty.
\eeao
Therefore $g(x) = o(a(x))$ in agreement with condition
\eqref{eq:4} which also implies that $g(x)/x\to 0$ since $S(x)\to\infty$
for $F\in\MDA(\Lambda)$ with infinite right endpoint.
 Moreover, we conclude from
\eqref{gumbel} that for any $c\in \R$,
\beam\label{eq:march16b}
\lim_{\xto} \dfrac{\P(X>x-c \,g(x))}{\P(X>x)}=1\,,
\eeam
i.e., $F$ is $(c g)$-insensitive for any $c\in \R$.
The latter condition will be frequently used in the remainder of this paper.
On the other hand, the first condition in \eqref{eq:3} ensures that
$g(x)$ increases not too slowly.

\ble\label{L1} If \eqref{eq:3} holds then
for any $\varepsilon>0$,
\beam\label{eq:2}
\lim_{\nto}\sup_{x>t_n} \dfrac{n \ov F(\varepsilon x) \ov F(\varepsilon g(x))}{\ov F(x)} = 0.
\eeam
\ele
\begin{proof}
Since $F\in {\bf LN}\cap \MDA(\Lambda)$ we have
$\lim_{\xto} S(x)/\log x = \infty.$ It follows from \eqref{eq:3} uniformly for $x>t_n$,
\beao
S(x)\sim S(g(x)) \ge S(g(t_n)) \ge S(\sqrt{n}) \gg \log n.
\eeao
Hence by slow variation of $S(x)$ and \eqref{eq:3}, uniformly for $x>t_n$,
\beao
 \dfrac{n \ov F(\varepsilon x) \ov F(\varepsilon g(x))}{\ov F(x)}
&=& \exp\big(\log n  - S(x)(1+o(1))\big)\to 0\,.
\eeao
\end{proof}
\bexam\label{remark2.5}\rm
If we choose $g(x)= x/S(x)$ and $S(x) = f(\log x)$
for a differentiable \regvary\ \fct\ $f$ with index $\alpha>1$  then \eqref{eq:4} and
\eqref{eq:3} are satisfied. Indeed, if $T(x)$ is slowly varying then,
according to  \cite{bojanic:seneta:1971},
the condition
\beam \label{bojanic}
\lim_{x\to\infty} \dfrac{xT^\prime(x)}{T(x)}\log T(x) = 0
\eeam
implies that for any $\rho\in \mathbb{R}$
\beam \label{seneta}
\lim_{x\to\infty} \dfrac{T(x \,T^{\rho}(x))}{T(x)} = 1
\eeam
holds and $T(x)=S(x)$ satisfies \eqref{bojanic}.
In particular, one can choose
$S_1(x) = c(\log x)^\alpha(1+o(1)),$ or  $S_2(x) = \exp(c(\log\log x)^\alpha)(1+o(1))$ for $c>0$,  and Lemma~\ref{L1} applies.
\par
The lognormal \ds\
is of type $S_1$ with $\alpha=2$; see \eqref{eq:feb20d}.
For this \ds\ we can choose $t_n\gg \sqrt{n}(\log n)^2.$
Therefore the conditions $g(t_n)/\sqrt{n}\to\infty$ and
\eqref{eq:feb20h} hold as well;
see the discussion in Section~\ref{subsec:ldsubexp}.
\eexam

\bre \label{reWeibull} Note that Weibull-type distributions do not satisfy conditions ${\bf C}_1$--${\bf C}_2$. Indeed, if a \ds\ $F$
has a tail $\ov F(x)=\exp(-x^{\alpha}L(x))$
for a \slvary\ \fct\ $L$ and $\alpha\in (0,1),$ then ${\bf C}_1$ implies that $g(x) \le C x^{1-\alpha}/L(x)$ as $\xto.$ Thus, $S(x)/S(g(x))\to\infty$ as $\xto$ and (\ref{eq:3}) is not satisfied.
\ere

\subsection{Time series models satisfying {\bf (C)}}
In the previous section we verified conditions ${\bf C}_1$--${\bf C}_2$
on some examples. These conditions depended only on the marginal \ds\
$F$ of $(X_i).$
In this section we provide some examples of \ts\ for which we can
verify condition $\bf C_3$ which depends on the pairwise dependence
structure of $(X_0,X_h)$, $h=1,\ldots,m$.
Here and in what follows,  $c$ denotes any positive constant whose
value is not of interest.

\bexam\label{exam:lognormal1}\rm
Let $\bfY = (Y_i)$ be a Gaussian $m$-dependent stationary sequence with mean $\mu,$ variance $\sigma^2>0$ and correlation function $\rho(h)<1$ for $h\neq 0.$ Consider a stationary sequence $\bfX = \ex^{b(\bfY)} = (\ex^{b(Y_i)}),$ where $b(x) = \sign (x)\,|x|^\alpha$,  $\alpha\in(0,2).$ We observe that for large $x$
\beao
S(x) = \dfrac{1}{2\sigma^2} (\log x)^{2/\alpha} (1+o(1)),
\eeao
thus $S(x)$ satisfies \eqref{bojanic} by Example 2.2 and then $g(x) = x/S(x)$ satisfies $\bf C_1.$ The conditions $g(t_n)/\sqrt{n}\to\infty$ and
\eqref{eq:feb20h} hold
with $t_n \gg \sqrt{n} (\log n)^{2/\alpha}.$ Indeed, according to \cite{rozovski:1993}, the large deviation result \eqref{eq:feb20h} holds with $t_n \gg \sqrt{n} (\log n)^{2/\alpha - 1}$ for $\alpha\in(0,1]$ and with $t_n \gg \sqrt{n} (\log n)^{1/\alpha}$ for $\alpha \in (1,2).$ Note also that for $\alpha=1$
the random vector $(X_1, \ldots, X_d),\,d\in \bbn,$ has  a
multivariate lognormal distribution in the sense of \cite{asmussen:rojas:2008}.
\par{
Next we verify $\bf C_3$. We assume $\mu=0$ and observe that
$\rho(h) = 0$ for $h> m$. An adapted version of Shibuya's classical estimate, \cite{shibuya:1959}, and
the tail-balance condition \eqref{tail-balance}
yield for $\vep>0$ and large $x$,
\beao
\P(|X_0| > \vep\, x\,, |X_h| > \vep\, g(x)) &\le&c\,
\P(X_0>\vep\,g(x), X_h>\vep\,g(x))\\ &=& c \,\P(\min(Y_0, Y_h) > (\log (\vep \,g(x)))^{1/\alpha})\\  &\le& c\, \P(Y_0 + Y_h > 2 \log(\vep \,g(x)))^{1/\alpha})\\
& = & c \ov \Phi\left(\dfrac{2(\log (\vep g(x)))^{1/\alpha} }{\sigma\sqrt{1+\rho(h)}}\right) = o\left(\ov \Phi\left(\dfrac{(\log x)^{1/\alpha} }{\sigma}\right)\right),
\eeao
where $\Phi$ is the standard normal distribution function.
In the last step we used the facts that $\rho(h)<1$ and
\beao
\dfrac{2(\log (\vep g(x)))^{1/\alpha} }{\sqrt{1+\rho(h)}}=
(\log x)^{1/\alpha}\,\dfrac{2 }{\sqrt{1+\rho(h)}}\,(1+o(1))\,.
\eeao
We conclude that for $h\ge 1$,
\beao
\P(|X_0| > \vep\, x\,, |X_h|> \vep\, g(x)) = o(\ov F(x)),\qquad \xto\,,
\eeao
and thus $\bf C_3$ is satisfied.}
\eexam

\bexam\rm
Let $(Y_i)$ be an iid sequence with common distribution given by
\beam \label{example2.4}
\P(Y>x) = \exp(-(\log x)^\alpha)\,,\qquad x > 1\,,
\eeam
for some $\alpha>1.$ The \seq\
\beao X_i = \min(a_0 Y_i, a_1 Y_{i+1}, \ldots, a_m Y_{i+m})
\eeao
for some positive $a_0,\ldots,a_m$ is $m$-dependent, stationary and has tail
\beao
\P(X>x) = \exp\big(-\sum_{i=0}^m S(x/a_i)\big)= \exp\big(-m\,S(x)(1+o(1))\big)\,.
\eeao
Thus the \ds\ of $X$ is also subexponential. This follows by checking
Pitman's condition, \cite{pitman:1980}: integrability
of the \fct\ $\exp(x F'(x)/\ov F(x)) F'(x)$  on $(0,\infty)$.
We verify
that {\bf (C)} holds with $g(x) = x/(\log x)^{\alpha}$. $\mathbf C_1$ is
immediate.
$\mathbf C_2$ follows by virtue of Example \ref{remark2.5}. It remains to
verify $\mathbf C_3$. Direct calculation yields for $\vep>0$ and $h=1,\ldots,m$,
\beao
&&\dfrac{\P(X_0>\varepsilon g(x)\,,X_h>\varepsilon x)}{\ov F(x)}\\ &\le&
\dfrac{\P\big( \min(a_0 Y_0,\ldots,a_{h-1}Y_{h-1})>\vep\,g(x)\,, \min(a_0Y_h,\ldots,
a_mY_{m+h})>\vep\,x\big)}{\P\big(\min (a_0Y_0,\ldots,a_mY_m)>x\big)}\\
&=& \exp\big( \sum_{i=0}^m S( x/a_i) - \sum_{i=0}^{h-1} S(\vep g(x)/a_i)
-\sum_{i=0}^m S(\vep x/a_i)
\big)\\
&=&\exp\big((1+o(1)) S(x) ((m+1)-(m+1+h))\big)\to 0\,,\qquad \xto\,.
\eeao
\eexam

\bexam\rm  Consider the stochastic volatility model
\beao
X_i = \sigma_i Y_i,
\eeao
where $(\sigma_i)$ is a stationary sequence with
$\P(a\le \sigma_1 \le b)=1,\, 0<a<1<b,$ and $(Y_i)$ is an iid
sequence with common distribution function $F_Y(x) = 1-\ex^{-S_Y(x)},$ \st\
$F_Y\in \MDA(\Lambda)\cap \mathcal{S}$, it satisfies
the tail-balance condition \eqref{tail-balance}, and  \eqref{bojanic} holds for $S_Y$. We also assume that the \ds\ $F$ of $X$ is subexponential.
This is not automatic even though it is easily verified that
$S(x)=S_Y(x)(1+o(1))$, hence $S(x)$ is \slvary , but this fact does not necessarily imply subexponentiality of $F$; see comments on p.~52 in
\cite{embrechts:klueppelberg:mikosch:1997}. Subexponentiality of $F$ can be
verified in simple situations, e.g. if $\sigma$ has a binomial \ds\ on $(a,b)$,
by using Pitman's aforementioned condition.
We choose as before $g(x) = x/S(x)$ and assume that it increases. Hence
$\mathbf C_1$--$\mathbf C_2$ are satisfied.
It remains to show $\mathbf C_3$. Applying the slow variation of $S(x)$,
the tail-balance condition \eqref{tail-balance} and $\mathbf C_2$,
we have for $h = 1, \ldots, m$ and $\varepsilon>0$,
\beao
\dfrac{\P(|X_0|>\varepsilon g(x)\,,|X_h|>\varepsilon x)}{\ov F(x)} & \le & \dfrac{\P(|Y_0|>\varepsilon g(x)/b\,,|Y_h|>\varepsilon x/b)}{\ov F_Y(x/a)}\\
& \le & c\dfrac{\ov F_Y(\varepsilon g(x)/b) \ov F_Y(\varepsilon x/b)}{\ov F_Y(x/a)} \\
& = & \exp(-S_Y(x) (1+o(1)))\to 0\,,\qquad \xto\,.
\eeao
\eexam
{
\subsection{Regularly varying stationary \seq s}
A random vector $\bfX$ with values in $\mathbb{R}^d$
and its distribution are {\it regularly varying with index $\alpha>0$}
if
\beao
\P\Big(\Big(\dfrac{\bfX}{x}\,,\dfrac{\bfX}{|\bfX|}\Big)\in\cdot\;\Big|\;|\bfX|>x\Big)\stw \P\big((Y,\bfTh) \in \cdot\big)\,,\qquad \xto\,,
\eeao
where $Y$ is Pareto distributed, $\P(Y>x)=x^{-\alpha}$, $x>1$, independent of
$\bfTh$;
see \cite{resnick:1987,resnick:2007} for some
reading on multivariate \regvar .
 \cite{davis:hsing:1995} introduced \regvary\ stationary \seq s
$(X_t)$ by assuming that each lagged vector $(X_0,\ldots,X_h)$, $h\ge 0$,
is \regvary\ with index $\alpha$. \cite{basrak:segers:2009}
characterized such \seq s by showing that \regvar\ of $(X_t)$ is equivalent
to the existence of a {\em spectral tail process} $(\Theta_t)$ defined via
the limit relations
\beam\label{eq:march17b}
\P(x^{-1}(X_{0},\ldots,X_h) \in \cdot\;\mid\; |X_0|>x) \stw
\P(Y\,(\Theta_{0},\ldots, \Theta_h) \in \cdot)\,,\qquad h\ge 0\,,\xto\,,
\eeam
where $Y$ is Pareto distributed and independent of $(\Theta_t)$. Obviously,
$|\Theta_0|=1$. If $\Theta_t=0$ a.s. for $t\ne 0$ then $(X_t)$ is called
\asy ally independent. In Section~\ref{subsec:regvar} we will work under
this assumption.
\par
We will work under the following set of conditions.
\subsection*{Condition (RV)}
\begin{enumerate}
\item[$\bf RV_1$]
The separating \seq\ $(t_n)$ satisfies
\beam\label{nlogn}
\lim_{\nto} \dfrac{t_n}{\sqrt{n\,\log n}} =\infty\,.
\eeam
\item[$\bf RV_2$] $F\in\RV(\alpha)$ for some $\alpha>2$ and satisfies the
tail-balance condition \eqref{tail-balance}.
\item[$\bf RV_3$] $(X_t)$ is $m$-dependent for some $m\ge 1$ and
satisfies
\beam\label{eq:march17a}
\lim_{\xto}\P(|X_h|>x\;\mid\; |X_0|> x)=0\,,\qquad
h=1,\ldots,m\,.
\eeam
\end{enumerate}
\par
Condition $\bf RV_2$ implies in particular that $\E[|X|^{2+\delta}]<\infty$
for $0<\delta<\alpha-2$. Moreover, $S(x)=-\log F(x)= \alpha \log x-\log L(x)$
for some \slvary\ \fct\ $L$. We conclude that any \fct\ $g$ satisfying
$g(x)/x\to 0$ as $\xto$ has the property
\beao
\lim_{x\to\infty}\dfrac{\ov F(x+ g(x))}{\ov F(x)} = 1\,.
\eeao
Condition $\bf RV_3$ implies the \asy\ independence of the \seq\ $(X_t)$, i.e.,
$\Theta_t=0$ a.s., $t\ne 0$, in \eqref{eq:march17b}. In particular, $(X_t)$
is \regvary\ with index $\alpha$. By \regvar\ we can rewrite \eqref{eq:march17a}
in the form
\beao
\lim_{\xto}\P(|X_h|>\vep \,x\;\mid\; |X_0|> \vep x)=0\,,\qquad
h=1,\ldots,m\,\,,\qquad \vep>0\,.
\eeao Condition $\bf RV_3$ is slightly stronger than the corresponding one
in \cite{mikosch:wintenberger:2016}
who proved their \ld\ result under the assumption that all
$\Theta_t$, $t=1,\ldots,m$, have an atom at zero. However, the proof
in this paper is direct in contrast to \cite{mikosch:wintenberger:2016}
who use techniques from the theory of \regvary\ processes.
In Section~\ref{subsec:ldsubexp} we mentioned that the best separating \seq\
in the iid \regvary\
case is $t_n=\sqrt{(\alpha-2)\,n\,\log n}$. Thus $\bf RV_1$ is
not too far away from the latter growth condition.
\bexam\rm
We consider the \sv\ model $X_t=\sigma_t\,Z_t$, $t\in\bbz$, where $(\sigma_t)$
is a positive stationary \seq\ independent of the iid \regvary\ \seq\ $(Z_t)$
with index $\alpha>0$. If $\E [\sigma^{\alpha+\delta}]<\infty$ for some $\delta>0$
then it is not difficult to see that $(X_t)$ is \regvary\  with index $\alpha$.
Moreover, it is \asy ally independent. Condition \eqref{eq:march17a}
can be verified as follows: for $h\ge 1$,
\beao
\P\big(|X_h|>x,\;|X_0|>x\big)& = &\P\big(\min (\sigma_h\,|Z_h|\,,\sigma_0\,|Z_0|)>x)\\
&\le & \P\big((\sigma_h\vee \sigma _0)\, (|Z_0|\wedge |Z_h|)>x\big)=:I(x)\,.
\eeao
We observe that $|Z_0|\wedge |Z_h|$ has \regvary\ tail with index $-2\alpha$.
By a result of  \cite{breiman:1965} we have
\beao
I(x)\sim \E\big[(\sigma_h\vee \sigma _0)^{2\alpha}\big]\,\P(|Z_0|\wedge |Z_h|>x)=\E\big[(\sigma_h\vee \sigma _0)^{2\alpha}\big]\,[\P(|Z|> x)]^2\,,
\eeao
provided $\E[\sigma^{2\alpha+\delta}]<\infty$ for some $\delta$.
The latter condition is satisfied e.g. if $\sigma$ has a lognormal \ds . This
is a standard assumption in financial \tsa ; see \cite{andersen:davis:kreiss:mikosch:2009}.
Since we also have $\P(|X|>x)\sim \E[\sigma^\alpha]\,\P(|Z|>x)$ relation
\eqref{eq:march17a} is immediate.
\eexam}
\section{Main results}\label{sec:main}\setcounter{equation}{0}
\subsection{$X$ has semi-exponential tails}
The following result is our main precise \ld\ result for
a stationary \seq\ with semi-exponential tails. The proof is given in Section~\ref{sec:proofmain}.
\bth \label{thm:main}
Consider an $m$-dependent stationary process $(X_i)$ with marginal \ds\ $F\in\MDA(\Lambda)\cap {\bf LN}$ for some $m\ge 1$.
Assume condition {\bf (C)}.
Then we have
\beam\label{main}
\lim_{\nto}\sup_{x>t_n\delta} \left|\dfrac{\P(S_n>x)}{n \ov F(x)} - 1\right| = 0, \quad \mbox{ for any}\ \delta>0.
\eeam
\ethe
An inspection of the proof of Theorem \ref{thm:main} shows
that it can be generalized in various directions.
Instead of $F\in\MDA(\Lambda)\cap {\bf LN}$ we may  require  $\E [|X|^{2+\delta}]<\infty$ for some $\delta>0,$ that $S(x)$ is slowly varying, $g$ satisfies ${\bf (C)}$ and
$F$ is $(\varepsilon g)$-insensitive for any $\varepsilon>0.$
Hence we may also take into consideration distributions with infinite moments,
in particular the class {\bf RV($\alpha$)} for some $\alpha>2$. However, the method
of proof does not allow one to get an ``almost''
optimal separating sequence $t_n \gg \sqrt{n \log n}$ under {\bf RV($\alpha$)}.
For this reason, we provide
Theorem \ref{thm:second}  under the  latter condition
which proves \eqref{main} for a best possible separating \seq .

\subsection{$X$ has \regvary\ tails}\label{subsec:regvar}
The following theorem complements the large deviation result
for $m$-dependent stationary regularly varying sequences
by \cite{mikosch:wintenberger:2016}. The methods
of proof are distinct and  do not make direct use of techniques for \regvary\
\seq s.
\bth\label{thm:second} Assume $(X_t)$ is an $m$-dependent stationary \seq\ which is \regvary\
with index $\alpha>2$ and condition $\bf (RV)$ is satisfied.
Then the large deviation result \eqref{main} holds.
\ethe
The proof of this result is given in Section~\ref{sec:proofthm32}.

\section{Linear process with subexponential noise}\label{subsec:linear}\setcounter{equation}{0}
Assume that $Z$ has a subexponential \ds\ $F_Z$ ($F_Z\in\mathcal S$)
in the sense that $Z_+$ has a subexponential \ds\  and a tail-balance condition holds:
\beam\label{eq:feb17a}
\dfrac{\P(Z>x)}{\P( |Z|>x)}\to p_{+}\,,\quad \dfrac{\P(-Z>x)}{\P( |Z|>x)}\to p_{-}\,, \qquad \xto\,,
\eeam
for some $p_+>0,$ $p_-\ge 0$ \st\ $p_++p_-=1$.
Throughout this section we assume $F_Z\in\MDA(\Lambda)\cap \mathcal S.$
Consider real coefficients $(\psi_j)$ \st\  $\psi_{j}=0$ for $j<0$,
$\max_j|\psi_j|=1$ and
\beam\label{eq:march15a}
\sum_{j=0}^\infty |\psi_j|^\delta<\infty \quad\mbox{ for some $\delta\in (0,1)$.}
\eeam
Let $k_{\pm}=\#\{j: \psi_j=\pm 1\}$,
$m_0=\sum_{j=0}^\infty\psi_j$ and $m_1=\sum_{j=1}^\infty|\psi_j|$ which are finite in view of \eqref{eq:march15a}. Then the infinite series
\beao
X=\sum_{j=0}^\infty \psi_{j} Z_{j}
\eeao
converges a.s. provided $(Z_t)$ is an iid \seq\ with generic element $Z$.
Indeed, $F_Z\in\MDA(\Lambda)\cap \mathcal S$ implies that $Z$ has finite first moment
and therefore $\E[\sum_{j=0}^\infty |\psi_j Z_j|]=(m_1+|\psi_0|)\E[|Z|]<\infty$.
\subsection{Tail behavior of $X$}
The following result was proved by  \cite{davis:resnick:1988}.
\ble\label{likeResnick}
If  \eqref{eq:march15a} and $F\in\MDA(\Lambda)\cap \mathcal S$ hold then
\beam\label{eq:feb11b}
\P(X>x)\sim k_+\,\P(Z>x)+ k_-\P(Z<-x)\sim (k_+p_+ + k_-p_-) \P(|Z|>x).
\eeam
\ele
We may conclude that the \ds\ of $X$ is tail-equivalent to $F_Z$. Hence it inherits subexponentiality.

\subsection{Large deviations of linear processes}
We consider the causal linear process
\beao
X_t= \sum_{j=0}^\infty\psi_jZ_{t-j}\,,\qquad t\in\bbz\,,
\eeao
with generic element $X$ for an iid \seq\ $(Z_t)$ with generic element $Z$,
$\E[Z]=0$ and
\beam\label{like-exp}
\sum_{i=1}^\infty \Big(\sum_{j=i}^\infty|\psi_j|\Big)^\delta <\infty\quad \mbox{for some $\delta\in (0,1)$}\,.
\eeam
This condition implies \eqref{eq:march15a}, and it is satisfied
if $\sum_{j=1}^\infty j\,|\psi_j|^\delta<\infty$.
Thus, by virtue of \eqref{eq:feb11b},  $X_t$ has a
subexponential \ds\ with tail balance condition.
The next result shows that a large deviation result for the iid subexponential
$(Z_t)$ with separating \seq\ $(t_n)$ implies a corresponding result with
separating \seq\ $(|m_0| t_n)$.
\par
In what follows, we write  $S_{n,Z}=Z_1+\cdots +Z_n$.

{
\bpr \label{pr4.2}
Consider a causal linear process $(X_t)$ with iid mean-zero
noise $(Z_t)$ with \ds\ $F_Z\in\mathcal S\cap \MDA(\Lambda),$ $m_0\ne 0$ and real
weights $(\psi_j)$ satisfying 
\eqref{like-exp}. Choose a function $g(x)$ such that $g(x) =  o(a(x))$
where $a(x)$ is an auxiliary function for $F_Z$ in the sense of \eqref{gumbel}.
\begin{enumerate}\item[\rm 1.]
Assume that for a separating \seq\ $(t_n)$ and a set $\Lambda_n\subset (|m_0|t_n(1+\delta),\infty)$ for any small $\delta>0$,
we have
\beam\label{eq:feb18a}
\limsup_{\nto}\sup_{x\in \Lambda_n}\dfrac{\P(m_1\,|Z|>g(x))}{n\,\P(|m_0||Z|>x)}=0\,.
\eeam
If $m_0>0$ we assume that for any small $\delta>0$,
\beam \label{large-dev}
\sup_{x>m_0 (1+\delta )t_n}\Big|\dfrac{\P(S_{n,Z}>x)}{n\,\P(Z>x)}- 1\Big|\to0\,,
\eeam
and if $m_0<0$ and $0<p_+<1$,
\beam\label{eq:feb19a}
\sup_{x>|m_0|t_n(1+\delta)}\Big|\dfrac{\P(-S_{n,Z}>x)}{n\,\P(Z\le -x)}- 1\Big|\to0\,.
\eeam
Then
\beam\label{eq:feb19b}
\limsup_{\nto}\sup_{x\in\Lambda_n}\Big| \dfrac{\P(S_n>x)}{n\,\P(|Z|>x/|m_0|)}-
p_{\pm}\1_{(0,\infty)}(\pm m_0)\Big|=0\,.
\eeam
\item[\rm 2.]
Assume $\psi_j=0$, $j>m$, for some $m\ge 1$,
and \eqref{large-dev} or
\eqref{eq:feb19a} hold according as $m_0>0$ or $m_0<0$. Moreover,
assume that there is a set $\Lambda_n$ \st\ $\Lambda_n\subset (|m_0|t_n(1+\delta),\infty)$ for any $\delta>0$ and, for $m_0'=\sum_{i=0}^m|\psi_i|$,
\beam\label{eq:feb19c}
\limsup_{\nto}\sup_{x\in\Lambda_n}\Big[
\dfrac{\P(m_0'|Z|>x)}{n\,\P(|m_0||Z|>x)}+\dfrac{[\P(m_0' |Z|>g(x))]^2}{\P(|m_0|\,|Z|>x)}\Big]
=0\,.
\eeam
Then \eqref{eq:feb19b} holds.
\end{enumerate}
\epr
\bre\label{remark4.3} We observe that $m_0'=|m_0|$ if either $\psi_j\ge 0$ for all $j$ or
$\psi_j\le 0$ for all $j$. In both situations, the first
ratio in
 \eqref{eq:feb19c} vanishes  for $\Lambda_n=(|m_0|t_n(1+\delta),\infty)$. The second ratio vanishes if $-2S(g(x)/|m_0|)+S(x/|m_0|)\to -\infty$ as $\xto$. This
condition holds if we can ensure that
$\sup_{x> t_n} |S(g(x))/S(x)-1|\to 0$.
The latter condition is satisfied for lognormal $Z$ if we choose
$g(x) = x/S(x)$ and
$t_n\gg \sqrt{n}\log n$.
\ere
\bexam\label{exam:4.4} \rm Condition \eqref{eq:feb18a} is quite restrictive. We illustrate
this for an iid \seq\ $(Z_i)$ with \ds\ given by
\beam\label{eq:march27a} \ov F_Z(x)=\P(Z> x) = \exp( - (\log x)^\alpha), \quad x>1,
\eeam
for $\alpha>1$ and $g(x) = \vep a(x),$ where
\beao \vep = \vep(x)\to 0 \mbox{ and } \vep(x) \log\log x\to \infty, \quad \xto.
\eeao
Calculation yields $a(x) \sim c x/(\log x)^{\alpha-1}$ for some $c>0$.
For convenience, we assume $m_0>0.$ We have
\beao
\dfrac{\P(m_1 Z > g(x))}{n\P(m_0 |Z|> x)} & = & \exp\big(-(\log(\vep a(x)/m_1))^\alpha + (\log(x/m_0))^\alpha - \log n\big)\\
& = & \exp\big( (\alpha-1) (\log x)^{\alpha-1} \log\log x (1 + o(1)) - \log n\big).
\eeao
For $\alpha\ge 2$ one can choose $t_n \gg \sqrt{n} (\log n)^{\alpha - 1}$ in
\eqref{eq:feb19a}; see the discussion in Example~\ref{exam:lognormal1}.
In this case $\Lambda_n$ is empty. For $\alpha \in (1,2)$ we can choose
\beao
\Lambda_n = (c_n, b_n), \quad c_n\gg \sqrt{n (\log n)^\alpha}, \quad b_n = \exp\Big( \Big(\frac{(1 - \delta)\log n}{(\alpha-1)\log\log n}\Big)^{1/(\alpha - 1)}\Big)
\eeao for arbitrarily small $\delta>0.$ In particular,
$c\,n\in \Lambda_n$ for any $c>0$.
\eexam
\bexam\rm We assume  $m_0'> m_0>0$ in \eqref{eq:feb19c}. In this case \eqref{eq:feb19c}
is as restrictive as  \eqref{eq:feb18a}. To illustrate this,
choose $F_Z$ as in \eqref{eq:march27a}.
As mentioned in Remark \ref{remark4.3}, the second summand in \eqref{eq:feb19c} vanishes for $x>t_n$ if
$\sup_{x>t_n} |S(g(x))/S(x) - 1| \to 0.$ We investigate
the first summand. We have
\beao
\dfrac{\P(m_0' Z> x)}{n\P(m_0 |Z|>x)} & =& \exp\big(- (\log x - \log m_0')^\alpha + (\log x - \log m_0)^\alpha - \log n\big)\\
 & = &  \exp \left(\alpha \log(m_0'/m_0) (\log x)^{\alpha - 1} (1+o(1)) - \log n \right).
\eeao
Thus we get similar restrictions as in Example~\ref{exam:4.4}.
The set $\Lambda_n$ is empty for $\alpha>2$. For $1<\alpha<2$
we can choose
\beao \Lambda_n = (c_n, b_n), \quad c_n\gg \sqrt{n (\log n)^\alpha}, \quad b_n = \exp\Big( \Big(\frac{(1-\delta)\log n}{\alpha \log(m_0'/m_0)}\Big)^{1/(\alpha - 1)}\Big)
\eeao
for arbitrarily small $\delta>0$,  and we observe that
$c\,n\in \Lambda_n$ for any $c>0.$ If $\alpha = 2,$ $\Lambda_n$ is not empty if $m_0'/m_0 < {\rm e}$ and contains the sequence $c\,n$ for any $c>0$
if $m_0\spp/m_0<\sqrt{{\rm e}}.$
\eexam

\begin{proof}[Proof of Proposition~\ref{pr4.2}] {\bf 1.}
We follow the ideas of the proof of Lemma A.5 in
 \cite{mikosch:samorodnitsky:2000}.
It will be convenient to write $\psi_j=0$ for $j\le 0$. We prove the result
for $m_0>0$; the case $m_0<0$ is analogous.
We start with the decomposition
\beao
S_{n}&=& \sum_{j=-\infty}^0 Z_j\beta_{n,j}+\sum_{j=1}^nZ_j\beta_{n,j}=:S_{n,1}+S_{n,2}\,,
\qquad\mbox{ where $
\beta_{n,j}= \sum_{i=1-j}^{n-j}\psi_i$.}
\eeao
We have
\beao
\P(S_{n,1}>x)&\le&  \P\Big(\sum_{j=0}^\infty|Z_j|\,\sum_{i=1+j}^{n+j}|\psi_i|>x\Big)\\
&\le& \P\Big(\sum_{j=0}^\infty|Z_j|\,\sum_{i=1+j}^{\infty}|\psi_i|>x\Big)
\le c\,\P(|Z|\,m_1>x)\,,
\eeao
where we used Lemma~\ref{likeResnick} in the last step. Indeed, the conditions
of this lemma are satisfied by virtue of \eqref{like-exp}.
We have
\beao
S_{n,2}=\sum_{j=1}^n Z_j\,\sum_{i=0}^{n-j}\psi_i\eqd \sum_{j=1}^n Z_j\,\sum_{i=0}^{j-1}\psi_j=m_0\,S_{n,Z}-\sum_{j=1}^n Z_j\,\sum_{i=j}^\infty\psi_j=
m_0\,S_{n,Z}-S_{n,21}\,.
\eeao
Applying  Lemma~\ref{likeResnick}, we obtain
\beao
\P(|S_{n,21}|>x)&\le &\P\Big(\sum_{j=1}^\infty |Z_j|
\sum_{i=j}^{\infty}|\psi_i|>x \Big)\le c\,\P(|Z|\,m_1>x)\,.
\eeao
By independence of $S_{n,1}$ and $S_{n,2}$ we observe that
\beao
\P(S_n>x)\le \P(S_{n,1}>x-g(x))+ \P(S_{n,2}>x-g(x)) +
\P(S_{n,1}>g(x))\,\P(S_{n,2}>g(x))\,.
\eeao
Hence for $m_0>0$, $x>t_n \,m_0(1+\delta)$ and sufficiently large $n$,
\beao
\P(S_n>x)
& \le &c\,\P(m_1\,|Z|>x-g(x))
+\P(S_{n,2}>x-g(x))\\
&&+c\,\P(S_{n,2}>g(x))\,\P(m_1\,|Z|>g(x))\\
& \le &c\,\P(m_1\,|Z|>x-g(x))+\P(m_0\,S_{n,Z}>x-2g(x))+ \P( -S_{n,21}>g(x))\\
&&+c\,\big(\P(m_0 S_{n,Z}>g(x)/2)+\P(-S_{n,21}>g(x)/2)\big)\P(m_1\,|Z|>g(x))\,.
\eeao
We conclude that
\beao\lefteqn{
\limsup_{\nto}\sup_{x\in\Lambda_n} \dfrac{\P(S_n>x)}{n\,\P(m_0|Z|>x)}}\\
&\le &\limsup_{\nto}\sup_{x\in\Lambda_n}\Big[
\dfrac{\P(m_0\,Z>x)}{\P(m_0 |Z|>x)} + c\,\dfrac{\P(m_1\,|Z|>g(x))}{n\,\P(m_0|Z|>x)} 
\Big]\\
&=&p_++ c\,\limsup_{\nto}\sup_{x\in\Lambda_n}\dfrac{\P(m_1\,|Z|>g(x))}{n\P(m_0|Z|>x)}=p_+\,.
\eeao
In the last steps we used \eqref{eq:feb18a}, the \ld\ result \eqref{large-dev} and the tail
balance condition \eqref{eq:feb17a}.
\par
We also have
\beao
\P(S_n>x)&\ge & \P(S_n>x\,,S_{n,1}\le g(x))
\ge \P(S_{n,2}>x+g(x)\,,S_{n,1}\le g(x))\\
&=&\P(S_{n,2}>x+g(x))\,\big(1-\P(S_{n,1}> g(x))=\P(S_{n,2}>x+g(x))(1+o(1))\,.
\eeao
Thus it suffices to find a lower  bound for
\beao
\P(S_{n,2}>x+g(x))&\ge &\P(m_0\,S_{n,Z}-S_{n,21}>x+g(x)\,,|S_{n,21}|\le  g(x))\\
&\ge &\P(m_0\,S_{n,Z}>x\,,|S_{n,21}|\le g(x))\\
&\ge &\P(m_0\,S_{n,Z}>x)- \P(|S_{n,21}|>  g(x))\\
&\ge & n\,\P(m_0\,Z>x)(1+o(1))-c\,\P(m_1 |Z|>g(x))\,.
\eeao
Therefore
\beao
\liminf_{\nto} \sup_{x\in\Lambda_n} \dfrac{\P(S_n>x)}{n\,\P(m_0|Z|>x)}
&\ge &\lim_{\nto} \sup_{x\in \Lambda_n}\dfrac{\P(m_0\,Z>x)}{\P(m_0|Z|>x)}\\&&
-c\,\limsup_{\nto} \sup_{x\in\Lambda_n}\dfrac{\P(m_1 |Z|>g(x))}{n\,\P(m_0|Z|>x)}=p_-\,.
\eeao
{\bf 2.}  We again assume $m_0>0$. In this case we have
\beao
S_n&=& \sum_{j=-m+1}^0Z_j \sum_{i=1-j}^m \psi_i + \sum_{j=n-m}^n Z_j\sum_{i=0}^{n-j}\psi_i+\sum_{j=1}^{n-m-1} Z_j m_0
=:T_{n,1}+T_{n,2}+T_{n,3}\,.
\eeao
Hence,
\beao
\P(S_n>x)&\le &\P(T_{n,1}+T_{n,2}>x-g(x))+ \P(T_{n,3}>x-g(x))\\&&+ \P(T_{n,1}+T_{n,2}>g(x))\P(T_{n,3}>g(x))\,.
\eeao
We have by Lemma~\ref{likeResnick} for sufficiently large $x$,
\beao
\P(T_{n,1}+T_{n,2}>x)&\le & \P\Big(\sum_{j=-m+1}^0|Z_j| \sum_{i=1-j}^m |\psi_i|+
\sum_{j=1}^{m+1}|Z_j| \sum_{i=0}^{j-1}|\psi_i|>x\Big)\\
&\le & c\,\P(m_0'\,|Z|>x)\,.
\eeao
Under \eqref{eq:feb19c},
\beao\lefteqn{
\limsup_{\nto}\sup_{x\in\Lambda_n}\dfrac{\P(S_n>x)}{n\,\P(m_0\,|Z|>x)}}\\
&\le & \limsup_{\nto}\sup_{x\in\Lambda_n} \Big[
\dfrac{\P(T_{n,1}+T_{n,2}>x-g(x))}{n\,\P(m_0|Z|>x)} + \dfrac{\P(T_{n,3}>x-g(x))}{n\,\P(m_0|Z|>x)}\\&&\hspace{5cm} + \dfrac{\P(T_{n,1}+T_{n,2}>g(x))\P(T_{n,3}>g(x))}{n\,\P(|m_0||Z|>x)}\Big]\\
&\le &\limsup_{\nto}\sup_{x\in\Lambda_n} \Big[\dfrac{\P(m_0' |Z|>x)}{n\,\P(m_0\,|Z|>x)}+\dfrac{\P(m_0\,Z>x)}{\P(m_0\,|Z|>x)}
+\dfrac{[\P(m_0' |Z|>g(x))]^2}{\P(m_0\,|Z|>x)}
\Big]=p_+\,.
\eeao
As regards the lower bound, we have uniformly for $x>m_0\,t_n\,(1+\delta)$,
\beao
\P(S_n>x)&\ge &\P(T_{n,3}>x+g(x))\,\P(T_{n,1}+T_{n,2}> -g(x))\\
&=&\P(T_{n,3}>x+g(x))\,(1-o(1))\\
 &\sim & n\,\P(m_0 \,Z>x+g(x))\sim p_+\,  n\,\P(m_0 \,|Z|>x)\,.
\eeao

\end{proof}}

\section{Application to a large sample covariance matrix}\setcounter{equation}{0}\label{sec:matrix}
Consider a real-valued field $(X_{it})$. We assume that
the rows $(X_{it})_{t\in\bbz}$, $i=1,2,\ldots,$ constitute iid stationary $m$-dependent \seq s. We observe the matrix $\bfX=(X_{it})_{i=1,\ldots,p;t=1,\ldots,n}$.
The corresponding sample covariance matrix is given by
\beao
\bfX\bfX^\top = \Big(\sum_{t=1}^nX_{it}X_{jt}\Big)_{i,j=1,\ldots,p}=: (S_{ij}^{(n)})_{i,j=1,\ldots,p}\,.
\eeao
We assume that $p=p_n\to\infty$. In what follows, $X$ stands for a generic
element of the field with \ds\ $F$, and we also write $(X_i)$ for an iid \seq\
with common \ds\ $F$.
\subsection{The case $F\in \RV(\alpha)$}
\ble\label{lem:march28a} Assume the following conditions:
\begin{itemize}
\item
$X\in{\RV}(\alpha)$ for some $\alpha>4$, in particular there is
$(c_n)$ \st\ $n\,\P({X^2}>c_n)\to 1$ and
$c_n^{-1}\max_{i=1,\ldots,n} X_i^2\std Y\sim \Phi_{\alpha/2}.$
\item The \asy\ tail relations are valid:
\beam
&&\P\big(S_{11}^{(n)}-n\,\E[X^2]>c_{np}x)\sim n\, \P(X^2>c_{np}x)\,,\qquad x>0\,,\label{eq:march29d}\\
&&\P\big(|S_{12}^{(n)}-n\,(\E[X])^2|> c_{np}x\big)\le c\, n\,\P(|X_{1}X_{2}|>c_{np}x)=o(p^{-2}),\, x>0.\label{eq:march29e}
\eeam
\end{itemize}
Then the following limit relations hold:
\beam
&&c_{np}^{-1} \max_{1\le i<j\le p} |S_{ij}^{(n)}-n\,(\E[X])^2|\stp 0\,,\label{eq:march29b}\\
&&c_{np}^{-1}\max_{i=1,\ldots,p} \big(S_{ii}^{(n)}-n\,\E[X^2]\big)\std Y.
\label{eq:march29a}
\eeam
\ele
\begin{proof}[Proof of Lemma~\ref{lem:march28a}]
By assumption \eqref{eq:march29d} we have for any $x>0$,
\beao
p\,\P(S_{11}^{(n)}-n\,\E[X^2]>c_{np}x)\sim (np)\,\P(X^2>c_{np}x)
\to x^{-\alpha/2}\,.
\eeao
The \rv s $(S_{ii}^{(n)})$ are iid and therefore \eqref{eq:march29a} holds \fif\
the latter relation does.
\par
Next we show that \eqref{eq:march29b} holds.
We have for any positive $x$,
\beao
\lefteqn{\P\Big(c_{np}^{-1} \max_{1\le i<j\le p} \big|(S_{ij}^{(n)}-n(\E[X])^2)
\big|>x
\Big)}\nonumber\\
&\le &\P\Big(\max_{1\le i<j\le p} (S_{ij}^{(n)}-n\,(\E[X])^2)>c_{np}x\Big)\nonumber\\
&&+\P\Big(\max_{1\le i<j\le p} (-S_{ij}^{(n)}+n\,(\E[X])^2)>c_{np}x\Big)=:I_1+I_2\,.
\eeao
We restrict ourselves to prove $I_1\to 0$. We have by assumption
\beao
I_1&\le & p^2\,\P\big(S_{12}^{(n)}-n(\E[X])^2>c_{np}\,x\big)
\le c\,p^2\,n\,\P(|X_{1}X_{2}|>c_{np}x\,)\to 0\,.
\eeao
\end{proof}
%

\bexam\rm
Assume that $X\in {\rm RV}(\alpha)$ for some $\alpha>4$
and the conditions  of Theorem~\ref{thm:second} are satisfied
for the \seq s $(X_{1t}^2)$ and $(X_{1t}X_{2t})$ with the same separating
\seq\ $(t_n)$ satisfying $t_n\gg \sqrt{n\,\log n}$.
Thus, $X^2$ is \regvary\ with
index $\alpha/2>2$ and $\P(X_1X_2>x)\sim x^{-\alpha}l(x)$ for some \slvary\ $l$;
see  \cite{embrechts:goldie:1980}.
We choose $(c_n)$ \st\ $n\,\P(X^2>c_n)\to 1$, i.e.,
$c_n=n^{2/\alpha}\ell(n)$ for some \slvary\ $\ell$.
We take $(p_n)$ such that $p= n^{\beta}$ with $\beta>\alpha/4 - 1,$ then we have $c_{np} \gg t_n.$
An application of Theorem~\ref{thm:second} yields for $x>0$,
\beao
\P\big(S_{11}^{(n)}-n\,\E[X^2]>c_{np}\,x\Big)\sim n\,\P(X^2>c_{np}\,x)\,.
\eeao
This is the desired relation \eqref{eq:march29d}. Next we consider
\beao
q_n=\P\big(|S_{12}^{(n)}-n\,(\E[X])^2|>c_{np}\,x\big)\,.
\eeao
By Theorem~\ref{thm:second} we have for some \slvary\ $\wt l$,
\beao
q_n\sim c\,n\,\P(|X_1X_2|>c_{np}\,x)= c\,n (np)^{-2} \wt l(np)=
n^{-1} p^{-2} \wt l(np)=o(p^{-2})
\eeao
provided  $\wt l(np)/n\to 0$. This condition is satisfied since we chose $p=n^\beta.$ Thus we have the desired relation
\eqref{eq:march29e}. We conclude that the limit relations
\eqref{eq:march29b} and \eqref{eq:march29a} for the maxima of the
diagonal and off-diagonal terms $S_{ii}^{(n)}$ and $S_{ij}^{(n)}$, $i\ne j$, hold.

\eexam

\subsection{The case $F\in\MDA(\Lambda)\cap \mathcal S$}
\ble\label{lem:auxy}
Assume the  following conditions:
\begin{itemize}
\item
The \ds\ of $X^2$ is in $\MDA(\Lambda)$, i.e., there exist constants $c_n>0$ and
$d_n\in\bbr$ \st\
$c_n^{-1}(\max_{i=1,\ldots,n} X_i^2-d_n)\std Y$ with standard Gumbel limit.
\item The \asy\ tail relations are valid:
\beam
\P\big(S_{11}^{(n)}-n\,\E[X^2]>c_{np}x+d_{np})\sim n\, \P(X^2>c_{np}x+d_{np})\,,\qquad x\in\bbr\,,\label{eq:feb12cond1}\\
\P\big(S_{12}^{(n)}-n\,(\E[X])^2> c_{np}x+d_{np}\big)\le c\, n\,\P(X_{1}X_{2}>c_{np}x+d_{np})=o(p^{-2}),\ x\in\bbr\,.\label{eq:feb12cond2}
\eeam
\end{itemize}
Then the following limit relations hold:
\beam\label{eq:feb12a}
&&c_{np}^{-1} \max_{1\le i<j\le p} ((S_{ij}^{(n)}-n\,(\E[X])^2)-d_{np})\stp -\infty\,,\qquad i\ne j\,,\\
&&c_{np}^{-1}\max_{i=1,\ldots,p} \big((S_{ii}^{(n)}-n\,\E[X^2])-d_{np}\big)\std Y\,. \label{eq:feb12b}
\eeam
\ele

\begin{proof}[Proof of Lemma~\ref{lem:auxy}]
By assumption \eqref{eq:feb12cond1} we have for any $x$,
\beao
p\,\P(S_{11}^{(n)}-n\,\E[X^2]>c_{np}x+d_{np})\sim (np)\,\P(X^2>c_{np}x+d_{np})
\to \ex^{-x}\,.
\eeao
The \rv s $(S_{ii}^{(n)})$ are iid and therefore \eqref{eq:feb12b} holds \fif\
the latter relation does.
\par
By assumption \eqref{eq:feb12cond2} we have for any $x\in\bbr$,
\beao
\P\Big(\max_{1\le i<j\le p} (S_{ij}^{(n)}-n\,(\E[X])^2)>c_{np}\,x+d_{np}
\Big)
&\le &p^2 \,\P(S_{12}^{(n)}-n\,(\E[X])^2>c_{np}\,x+d_{np})\\
& \le & c\,p^2\,n\,\P(X_1X_2>c_{np}\,x+d_{np})\to 0\,.
\eeao
Relation \eqref{eq:feb12a} follows.
\end{proof}
\bexam\rm Assume that $(X_{1t})$ is $m$-dependent stationary with a lognormal
generic element $X$ and the conditions of
Example~\ref{exam:lognormal1} are met for $\alpha=1$. We standardize the
marginal \ds\ \st\ $X\eqd \ex^N$ for a standard normal \rv\ $N$. Thus,
$X^2\eqd \ex^{2N}$ and $X_1X_2\eqd \ex^{\sqrt{2}N}$ for independent copies $(X_i)$ of $X$.
According to Example~\ref{exam:lognormal1} we can apply
Theorem~\ref{thm:main} to both $(X_{it}^2)$ and $(X_{it}X_{jt})$ for $i\ne j$
and in both cases we can choose any separating \seq\
$t_n\gg \sqrt{n}(\log n)^2$.
\par
We set
\beao
c_{n}=2\,(2\log n)^{-1/2}d_n\,,\qquad d_n=\exp\big(2\big(\sqrt{2\log n}-(\log (4\pi)+\log\log n)/(2\sqrt{2\log n})\big)\big)\,.
\eeao
It is well known that $n\,\P(X^2>c_n x+d_n)\to \ex^{-x}$ for any $x\in\bbr$;
see \cite{embrechts:klueppelberg:mikosch:1997}, Example 3.3.31.
We take $(p_n)$ \st\  $p\gg n^{-1}\exp(C\,(\log n)^2)$ for some $C>1/16$.
Since $c_n= o(d_n)$ we have
\beao
c_{np}x+d_{np}\gg t_n\quad\mbox{ for any negative $x$.}
\eeao
Therefore,
\eqref{eq:feb12cond1} follows from Theorem~\ref{thm:main}.
\par
Next we verify \eqref{eq:feb12cond2}. To get the first bound in this
relation we apply
Theorem~\ref{thm:main}. Again observing that $c_n=o(d_n)$, we have
for $x\in\bbr$,
\beao
\P\big(S_{12}^{(n)}-n\,(\E[X])^2>c_{np}x+d_{np}\big)&\sim&
n\, \P(X_1X_2> c_{np}x+d_{np})\\
&=&p^2\,n\,\P(N>\log (c_{np}x+d_{np})/\sqrt{2})\\
& = & p^2\,n\,\P\big(N>2\sqrt{\log (np)}(1+o(1))\big)\\
&\sim &p^2\,n\,\dfrac{\ex^{-2\log (np) (1+o(1))}}{\sqrt{\pi} \log d_{np}}\to 0\,.
\eeao
Therefore \eqref{eq:feb12cond2} holds. We conclude that
the statements of Lemma~\ref{lem:auxy} are valid.
\eexam

\section{Proof of Theorem~\ref{thm:main}}\setcounter{equation}{0}\label{sec:proofmain}
For later use we recall two classical inequalities.
Consider a \seq\ $(X_i)$ of independent mean-zero \rv s, set
$S_n = \sum_{i=1}^n X_i$ and $\sigma^2_n = \var(S_n).$
\begin{itemize}
\item
{\bf Prokhorov's inequality}  (\cite{prokhorov:1959}; see \cite{petrov:1995}, p.~77)
If $|X_i| \le c$ a.s. for $i=1,\ldots, n$
and some constant $c$ then
\beao
\P(S_n>x) \le \exp\left( - \dfrac{x}{2c}{\rm arsinh}\big(\dfrac{cx}{2\sigma^2_n}\big)\right), \qquad x>0,
\eeao
where ${\rm arsinh}(y) = \log\big(y + \sqrt{y^2+1}\big).$
\item {\bf Fuk-Nagaev's inequality} (\cite{fuk:nagaev:1971,fuk:nagaev:1976}; see \cite{petrov:1995}, p.~78)
If $\E[|X_i|^p]<\infty$ for some $p\ge 2$, $i=1,\ldots,n$, $m_{p,n}= \sum_{i=1}^n
\E[|X_i|^p]$, then for constants $c_p,d_p>0$ only depending on $p$,
\beao
\P(S_n>x) \le c_p\,m_{p,n}\,x^{-p}+\ex^{-d_p\,(x/\sigma_n)^2}\,,\qquad x>0\,.
\eeao
\end{itemize}
\subsection*{The lower bound}
We have
\beao
\{S_n>x\}\supset \bigcup_{i=1}^n \{|S_n-X_i|\le g(x)\,,X_i>x+ g(x),
|X_j|\le g(x), 1\le j\ne i\le n\}\,.
\eeao
The events on the \rhs\ are disjoint. Therefore
\beam
\P(S_n>x)&\ge& \sum_{i=1}^n \P\Big(|S_n-X_i|\le g(x)\,,X_i> x+ g(x),
\max_{j\ne i}|X_j|\le g(x)\Big)\nonumber\\
&=&\sum_{i=1}^n \P\Big(X_i>x+ g(x),
\max_{j\ne i}|X_j|\le g(x)\Big)\nonumber\\
&&-\sum_{i=1}^n\P\Big(|S_n-X_i|>  g(x)\,,X_i>x+ g(x),
\max_{j\ne i}|X_j|\le  g(x)\Big)\nonumber\\
&=&\sum_{i=1}^n \P(X_i>x+ g(x))-\sum_{i=1}^n
\P\Big(X_i>x+ g(x)\,,\max_{j\ne i}|X_j|> g(x)\Big)\\
&&-\sum_{i=1}^n\P\Big(|S_n-X_i|> g(x)\,,X_i>x+g(x),
\max_{j\ne i}|X_j|\le g(x)\Big)\nonumber\\
&=&J_1(x)-J_2(x)-J_3(x)\,.\label{eq:march17e}
\eeam
We have
$\sup_{x>t_n}\big|J_1(x)/(n\,\ov F(x))-1\big|\to 0$ as $\nto$ and
\beao
\sup_{x>t_n}\dfrac{J_2(x)}{n\ov F(x)}&\le & \sup_{x>t_n}\sum_{i=1}^n \dfrac{\P\Big(X_i>x+ g(x), \max_{j\ne i:|j-i|\le m} |X_j|>g(x)\Big)}{n\ov F(x)}\\&&+ \sup_{x>t_n}\sum_{i=1}^n \dfrac{\P(X>x+ g(x))}{n\,\ov F(x)} \P\Big(\max_{j\ne i:|j-i|> m} |X_j|> g(x)\Big)\\
&\le & 2\,\sup_{x>t_n}\sum_{h=1}^m \dfrac{\P(X_0>x+ g(x)\,,|X_h|> g(x))}{\ov F(x)}\\
&&+\sup_{x>t_n}\dfrac{\P(X>x+ g(x))}{\ov F(x)} \P\Big(\max_{i=1,\ldots,n}|X_i|>
 g(x)\Big)\\
&\leq &o(1)+ \sup_{x>t_n}n\,\P(|X|>g(x))=o(1)\,,
\eeao
where the latter relation follows by $\bf C_3$ 
and since
$g(t_n)/\sqrt{n}\to\infty$ holds. Thus, it is enough to show that $J_3(x) \to 0$ as $\nto$ to derive the required lower bound for $\P(S_n>x).$
\par
In the sequel, we will use the notation,
\beam\label{eq:march16a}
\wh X_j= X_j\1(|X_j|\le g(x))\,,\qquad  \wh S_n^{(i)}= \sum_{j=1,j\ne i}^n \wh X_j\,.
\eeam
Hence,
\beao
\nonumber
J_3(x)&=&\sum_{i=1}^n \P\big(|\wh S_n^{(i)}|> g(x)\,,X_i>x+ g(x)\big)\\
\nonumber
&\le &\sum_{i=1}^n \P\Big( \Big|\sum_{1\le t\ne i\le n\,,|t-i|\le m} \wh X_t  \Big|> g(x)/2\,, X_i>x+ g(x)\Big)\\
&&+\sum_{i=1}^n \P\Big( \Big|\sum_{1\le t\le n\,,|t-i|>m} \wh X_t  \Big|> g(x)/2
\Big)\, \P(X>x+ g(x))\\
\nonumber
&=&J_{31}(x)+J_{32}(x)\,.
\eeao
We have by $\bf C_3$\,,
\beao
\dfrac{J_{31}(x)}{n\,\ov F(x)}&\le & 2\sum_{h=1}^m \dfrac{\P(|X_h|> g(x)/(4m)\,, X_0>x+ g(x)) }{\ov F(x)}=o(1)\,.
\eeao
 Finally, we deal with $J_{32}(x).$ Since $X$ has mean zero, we derive for arbitrary $\delta>0$
\beao
n\,|\E[\wh X]| = n\,|\E [X \1(|X|<g(x))]| = n\,|-\E [X \1(|X|>g(x))]| \leq \frac{n\E[ |X|^{2+\delta}]}{g^{1+\delta}(x)}.
\eeao
We deduce from \eqref{eq:3} and the fact that $\E [|X|^{2+\delta}]<\infty$,
\beam \label{mean_markov}
n|\E[\wh X]| \le n\dfrac{\E [|X|^{2+\delta}]}{g^{1+\delta}(t_n)} = o\big(g^{1-\delta}(t_n)\big), \quad \nto.
\eeam
Write $\wt g_r(x) = g(x)/(2m)- \#N_r\,\E[\wh X]$ where
\beao
N_r=\{1\le t\le n:t\equiv r({\rm mod}\;m), |t-i|>m\}\,,
\eeao
and observe that $n\, |\E[\wh X]| = o(g^{1-\delta}(x))$ and
$|\wh X-\E[\wh X]|\le 2g(x)$.
Using the $m$-dependence and Prokhorov's inequality,
we have for iid copies $(X_i\spp)$ of $X$ and large $n$,
\beao
\frac{J_{32}(x)}{n\ov F(x)} &\le & \sum_{i=1}^n\sum_{r=1}^m\P\Big(
\Big|\sum_{t\in N_r} \wh X_t\spp\Big|> g(x)/(2m)  \Big)
\dfrac{\P(X>x+ g(x))}{n \ov F(x)}\\
&\le & \sum_{i=1}^n\sum_{r=1}^m \P\Big(\Big|\sum_{t\in N_r} (\wh X_t\spp -\E [\wh X])\Big| > \wt g_r(x) \Big)\,\dfrac{\P(X>x+ g(x))}{n \ov F(x)}\\
& \le & c \sum_{i=1}^n\sum_{r=1}^m\exp \left( - \frac{\wt g_r(x)}{4g(x)} {\rm arsinh}\left(\frac{ 2 g(x)\wt {g}_r(x)}{2\# N_r \var(\wh X)}\right)\right)
\,\dfrac{\P(X>x+ g(x))}{n \ov F(x)}
\\
&\le & c\,\, m \exp \Big(-\Big(\dfrac 1{8m} +o(1)\Big)\, \log \Big((1+o(1))\frac{ m\,(g(x))^2}{2\,n \var (X)}\Big) \Big) \to 0\,.
\eeao
In the last step we used that $(g(x))^2/n\ge (g(t_n))^2/n\to \infty$; see
\eqref{eq:3}.
\subsection*{The upper bound}
Consider the following disjoint partition of $\Omega$: 
\beao
B_1&=& \bigcup_{1\le i<j\le n} \{|X_i|>g(x),|X_j|>g(x)\}\,,\\
B_2&=& \bigcup_{i=1}^n \{|X_i|>g(x)\,,\max_{j=1,\ldots,n,i\ne j}|X_j|\le g(x)\}\,,\\
B_3&=&\Big\{\max_{j=1,\ldots,n}|X_j|\le g(x)\Big\}\,.
\eeao
\subsection*{The bound on $B_1$}
We observe that for any $\xi\in(0,1)$,
\beao
\P\big(\{S_n>x\}\cap B_1\big)
&\le & \sum_{1\le i<j\le n} \P\big(S_n>x, |X_i|>g(x)\,,|X_j|>\,g(x)\big)\\
& \le & \sum_{1\le i<j\le n} \P\big(S_{ij}^{(1)} >  \xi x, |X_i|>g(x)\,,|X_j|>\,g(x)\big) \\
&& + \sum_{1\le i<j\le n} \P\big( S_{ij}^{(2)} > (1-\xi) x, |X_i|>g(x)\,,|X_j|>\,g(x)\big)\\
& =: & R_1(x) + R_2(x),
\eeao
where
\beao S_{ij}^{(1)} &=& \sum_{h\le n:\,|i-h|\wedge |j-h|> m} X_h = \sum_{r=1}^{m} (S_{ij}^{(r)})\spp
, \quad S_{ij}^{(2)}=\sum_{h\le n:\,|i-h|\wedge |j-h|\le  m} X_h\,,\\
(S_{ij}^{(r)})\spp&=& \mbox{$\sum_{h\in Q_{ij}^{(r)}} X_h$}\,,\qquad Q_{ij}^{(r)}=\{h\le n:|i-h|\wedge|j-h|> m, h \equiv r({\rm mod}\,m)\}.
\eeao
For a given $r$, the summands in $(S_{ij}^{(r)})\spp$ are independent due to
 $m$-dependence and also independent of $X_i$, $X_j.$ We have
$\# Q_{ij}^{(r)}\le n/m$ while the number of summands in $S_{ij}^{(2)}$ does not exceed $4m+2.$
Thus by the large deviation result \eqref{eq:feb20h}, $m$-dependence and stationarity,
\beao
\dfrac{R_1(x)}{n\,\ov F(x)} &\le & \sum_{1\le i<j\le n} \sum_{r=1}^m\dfrac{
\P\big((S_{ij}^{(r)})\spp >  \xi x/m\big)\,\P\big( |X_i|>g(x)\,,|X_j|>\,g(x)\big)}{n\,\ov F(x)}\\
& \sim & \dfrac{\ov F(\xi x/m)}{\ov F(x)}\sum_{1\le i<j\le n}\P\big(|X_i|>g(x)\,,|X_j|>\,g(x)\big)\\
&=&\dfrac{\ov F(\xi x/m)}{\ov F(x)}\sum_{h=1}^{n-1} (n-h)\P\big(|X_0|>g(x)\,,|X_h|>\,g(x)\big)\\
&\le & n \dfrac{\ov F(\xi x/m)}{\ov F(x)}\,\sum_{h=1}^{m}\P\big(|X_0|>g(x)\,,|X_h|>\,g(x)\big)+\dfrac{\ov F(\xi x/m)}{\ov F(x)}\, [n\,\P(|X|>g(x))]^2\\
&=:&R_{11}(x)+R_{12}(x)\,.
\eeao{
Applying the tail balance condition, $\bf C_3$ and \eqref{eq:2}, we have
\beao
\sup_{x>t_n}R_{11}(x)
\le  c\,m \sup_{x>t_n}\dfrac{n \ov F(\xi x/m) \ov F(g(x))}{\ov F(x)}\to 0\,.
\eeao
Since $g(t_n)/\sqrt{n}\to \infty$ and $\E[X^2]<\infty$ we also have
\beam \label{nPg}
\sup_{x>t_n}n\,\P(|X|>g(x)) \le \sup_{x>t_n}n\P(|X|>\sqrt{n}) \to 0\,.
\eeam
Hence, the tail balance condition, Lemma \ref{L1} and \eqref{nPg} immediately imply that $R_{12}\to 0.$}
\par
We have
\beao
R_2(x) &\le& \sum_{i=1}^{n-1} \Big(\sum_{j=i+1}^{(i+2m)\wedge n}+
\sum_{j-i>2m,j\le n}\Big)
 \P\big(S_{ij}^{(2)} > (1-\xi) x, |X_i|>g(x)\,,|X_j|>\,g(x)\big)\\
&=&R_{21}(x)+R_{22}(x)\,.
\eeao
We restrict ourselves to the study of $R_{21}(x)$; $R_{22}(x)$ can be treated
by similar methods. We note that $S_{ij}^{(2)}$ has \rep\
\beao
S_{ij}^{(2)} = \sum_{h=(i-m)\vee 1}^{(j+m)\wedge n} X_h.
\eeao
Observe that the number of summands in $S_{ij}^{(2)}$ does not exceed $4m+2$.
Therefore and by stationarity, taking care of the
cases $h = i$ and $h = j,$
\beam \label{eq:aug22}
R_{21}(x) & \le & \sum_{i=1}^{n-1} \sum_{j=i+1}^{(i+2m)\wedge n} \sum_{h=(i-m)\vee 1}^{(j+m)\wedge n}
\P\Big(|X_h|>\dfrac{(1-\xi)x}{ 4m+2}\,,
|X_i|>g(x),|X_j|>g(x)\Big)\\
& \le & c \,n\, \sum_{h=1}^m \P\Big(|X_0|>\dfrac{(1-\xi)x}{4m+2}\,,
|X_h|>g(x)\Big) + c\, n\, \ov F\Big(\dfrac{(1-\xi)x}{4m+2}\Big)\ov F(g(x))\,.
\eeam
By $\bf C_3$ and \eqref{eq:2} we conclude that
\beao
\limsup_{\nto} \sup_{x>t_n} \dfrac{R_{21}(x)}{n\ov F(x)} =0\,.
\eeao
Combining the previous bounds, we conclude that
\beao
\lim_{\nto}\sup_{x>t_n}\dfrac{\P\big(\{S_n>x\}\cap B_1\big)}{n\ov F(x)} =0\,.
\eeao
\subsection*{The bound on $B_2$}
\par{
Next we bound $\P(\{S_n>x\} \cap B_2)$.
Recall the notation $\wh X_j$ and $\wh S_n^{(i)}$ from \eqref{eq:march16a}.
Fix $b\in(0,1).$
Since $g(x)/x\to 0$ as $\xto$ we have
\beao
\P\big(\{S_n>x\}\cap B_2\big) & \le & \sum_{i=1}^n \P\big(X_i +
\wh S_n^{(i)} > x\,,|X_i|>g(x)\big) \\
 & = &  \sum_{i=1}^n\P\big(X_i + \wh S_n^{(i)}
> x\,,|X_i| \in (g(x), x - bx]\big) \\
&&+ \sum_{i=1}^n\P\big(X_i + \wh S_n^{(i)} > x\,,|X_i| \in (x - bx, x - g(x)]\big)\\ &&+ \sum_{i=1}^n\P\big(X_i + \wh S_n^{(i)} > x\,,|X_i| > x - g(x)\big)\\
&=:& I_1(x) +I_2(x) + I_3(x)\,.
\eeao
\subsubsection*{Bounding $I_1(x)$}
We show that $I_1(x) = o(n \ov F(x)).$
Similarly to the bound for $R_1(x)$ estimation, using the $m$-dependence,
we split $\wh S_n^{(i)}$  into $m$ sums of iid summands:
\beao
 \wh S_n^{(i)} = \sum_{r = 1}^{m} \wh S_{i,r}, \quad \wh S_{i, r} = \sum_{h\in Q^\ast_{i, r}} \wh X_h, \quad \mbox{ where }\quad Q^\ast_{i, r} = \{h\le n: h \equiv r ({\rm mod}\ m), h\neq i\}.
\eeao
In view of \eqref{mean_markov} we have $n|\E[\wh X]|\le g(x)$ for large $x$.
Moreover,
$|\wh X - \E[\wh X]|\le 2\,g(x)$  and $\# Q_{i, r} \le n/m.$
An application of Prokhorov's inequality for large $n$ yields}
\beao I_1(x) &\le & \sum_{i=1}^{n} \P\big( \wh S_n^{(i)} > b\,x\big)
\le \sum_{i=1}^{n} \sum_{r=1}^{m}\P\big(\wh S_{i, r} > bx/m\big) \nonumber \\
& \le& \sum_{i=1}^{n} \sum_{r=1}^{m}\P\big(\wh S_{i, r} - \# Q_{i, r} \E[\wh X] >
bx/m-g(x)\big) \nonumber \\
 & \le & \sum_{i=1}^{n} \sum_{r=1}^{m} \exp\Big(-\dfrac{bx/m-g(x)}{4 g(x)}
{\rm arsinh} \Big(\dfrac{2 g(x) (b x/m-g(x))}{2 \# Q_{i,r} \var(\wh X)}\Big)\Big)\\
 &\le& n\,m\exp\Big(- c\dfrac{x}{g(x)}
\log \Big(\dfrac{x\,g(x)}{ n}\Big)\Big)=o(n\,\ov F(x))
\eeao
uniformly for $x>t_n$. In the last step we used the bounds on $g(x)$ in
$\bf C_1$ and $\bf C_2$.
\subsubsection*{Bounding $I_2(x)$}
Write
\beao
\wt S_n^{(i)} = \sum_{t\notin[i-m,i+m]}\wh X_t.
\eeao
We observe that $|\wh X|\le g(x)$ and conclude by independence between
$X_i$ and $\wt S_n^{(i)}$, the tail-balance condition and integration by parts
that
\beao
I_2(x) &\le& \sum_{i=1}^{n} \P\Big(X_i+\wt S_n^{(i)}> x - 2\,m \,g(x)\,, |X_i| \in (x - bx, x - g(x)] \Big)\\
&\le& c\sum_{i=1}^{n}\int_{x-bx}^{x - g(x)} \P(\wt S_n^{(i)}> x - 2\,m\,g(x) - y)\,
d F(y)\\
&\le& c\,{\ov F(x - bx)}\sum_{i=1}^{n}\P(\wt S_n^{(i)}> bx - 2mg(x)) + c\int_{g(x)}^{bx}\sum_{i=1}^{n}\ov F(x-y)\,\P(\wt S_n^{(i)} \in dy) \\
&=:& c(I_{21}(x) + I_{22}(x))\,.
\eeao
In view of \eqref{eq:march16b}, for every $\delta>0$ there is $u_\delta$
\st
\beao
\sup_{y\ge u_\delta/(1-b)} \dfrac{\ov F(y)}{\ov F(y + g(y))} \le  \ex^{\delta}.
\eeao{
By a telescoping argument for sufficiently large $x$,
\beao
\dfrac{\ov F(x - bx)}{\ov F(x)} \le  \prod_{h=1}^{[bx/g(x)]} \dfrac{\ov F(x - hg(x))}{\ov F(x - (h-1)g(x))} \dfrac{\ov F(x-bx)}{\ov F(x- g(x)[bx/g(x)])}\le \ex^{\delta ([(bx)/g(x)]+1)}\,.
\eeao
Now the same argument as for $I_1(x)$ combined with $\bf C_1,$ \eqref{eq:3} and Prokhorov's inequality yields uniformly for $x>t_n$,
\beao
\dfrac{I_{21}(x)}{n\,\ov F(x)} &\le& \ex^{\delta ([bx/g(x)]+1)} \,m \,
\exp\left( -\dfrac{bx/m-3g(x)}{4 \,g(x)} \log\Big(\dfrac{m g(x)(bx/m-3g(x))}{ n\,\var(X)}\Big) \right) \\
&\le& \exp\left(-c\,\dfrac{x}{g(x)} \log\Big(\dfrac{x g(x)}{n}\Big)\right)\to 0\,,\qquad \nto\,.\eeao
A similar argument yields
\beao
\dfrac{I_{22}(x)}{n\,\ov F(x)} &\le & \dfrac{1}{n\ov F(x)}
\sum_{k=1}^{[bx/g(x)]} \int_{g(x)k}^{g(x)(k+1)}
\sum_{i=1}^{n}\ov F(x-y)\,\P(\wt S_n^{(i)} \in dy)\\
&\le &\sum_{k=1}^{[bx/g(x)]}  \dfrac{\ov F(x - (k+1)g(x))}{n\,\ov F(x)}\,\sum_{i=1}^{n}\P(\wt S_n^{(i)} \in g(x)\,(k, k+1])\\
&\le& \dfrac 1 n\sum_{k=1}^{\infty} \ex^{(k+1)\,\delta}\sum_{i=1}^{n}\P(\wt S_n^{(i)} > k\, g(x))\,\\
&\le& \sum_{k=1}^{\infty} \exp\left(- c k \log\Big(\dfrac{k g^2(x)}{n}\Big) + (k+1)\,\delta\right)\to 0\,,\qquad \nto\,.
\eeao
\subsubsection*{Bounding $I_3(x)$}
We have
\beao
\limsup_{x>t_n} \dfrac{I_3(x)}{n\,\ov F(x))} &\le & \limsup_{\nto}\sup_{x>t_n}\dfrac{1}{n\,\ov F(x)}\Big(\sum_{i=1}^n\P\big(X_i +\wh S_n^{(i)} > x\,,X_i > x - g(x)\big)\\&& + \sum_{i=1}^n\P\big(X_i + S_n^{(i)} > x\,,X_i < - x + g(x)\big)\Big) \\
&\le&  \limsup_{\nto}\sup_{x>t_n}\dfrac{\ov F(x - g(x))}{\ov F(x)} + \limsup_{\nto}\sup_{x>t_n} \dfrac{1}{n\,\ov F(x)}\sum_{i=1}^n \P\big(\wh S_n^{(i)} > 2x - g(x)\big)\\
&\le & 1+\lim_{\nto}\sup_{x>t_n}\exp\left(-c\dfrac{x}{g(x)}\log\Big(\dfrac{xg(x)}{n}\Big)\right)=1\,.
\eeao
The second term is bounded in the same was as $I_1$,
by exploiting the $m$-dependence, $\bf C_1$, \eqref{eq:3} and Prokhorov's inequality.
\par
Collecting the bounds for all $I_i(x)$, we obtain the desired relation
\beao
\lim_{\nto} \sup_{x>t_n}\dfrac{\P\big(\{S_n>x\}\cap B_2\big)}{ n \,\ov F(x)}\le 1\,.
\eeao
\subsection*{The bound on $B_3$}
It remains to show that $\P\big(\{S_n>x\}\cap B_3)=o(n\ov F(x))$.
We observe that
$\{S_n>x\}\cap B_3=\{\wh S_n>x\}$ where $\wh S_n=\sum_{i=1}^n\wh X_i$
and $|\wh X_i|\le g(x)$.
Now the same techniques as for bounding $I_1(x)$ apply.
We omit further details.
This finishes the proof of the upper bound. \qed
}
\section{Proof of Theorem~\ref{thm:second}}\setcounter{equation}{0}\label{sec:proofthm32}{
The proof is similar to the one of Theorem \ref{thm:main}. We follow the lines
of this proof and also use the same notation. We set $g(x)=g_\vep(x)=\vep\,x$
for any $\vep >0$.
\subsection*{The lower bound}
We start with the bound \eqref{eq:march17e}: $\P(S_n>x)\ge J_1(x)-J_2(x)-J_3(x)$.
For the first term we have
\beao
\lim_{\nto} \sup_{x>t_n} \dfrac{J_1(x)}{n\,\ov F(x)}&=&
\lim_{\nto} \sup_{x>t_n} \dfrac{n\,\ov F(x(1+\vep))}{n\,\ov F(x)}\\
&=& (1+\vep)^{-\alpha}\,,
\eeao
and the \rhs\ converges to 1 as $\vep\downarrow 0$. By \regvar\ the bound on $J_2(x)$
turns into
\beao
\sup_{x>t_n} \dfrac{J_2(x)}{n\,\ov F(x)}&\le & 2\sup_{x>t_n}\sum_{h=1}^m   \dfrac{
\P(X_0>(1+\vep)\,x\,,|X_h|>\vep x)}{\ov F(x)}\\
&&+\sup_{x>t_n}  \dfrac{\P(X>(1+\vep)\,x)}{\ov F(x)} \,
\P\Big(\max_{i=1,\ldots,n} |X_i|>\vep x\Big)\\
&\le &c\,\sup_{x>t_n}\sum_{h=1}^m
\P(|X_h|>\vep \,x\;\mid\;  |X_0|>\vep x)+
c\,n\,\P(|X|>t_n)\to 0\,.
\eeao
 Here we used condition \eqref{eq:march17a} for the
first term
and the facts that
$t_n\gg \sqrt{n}$ and $\var(X)<\infty$ for the second term.
Next we consider the bound $J_3(x)\le J_{31}(x)+J_{32}(x)$.
The relation $\lim_{\nto}\sup_{x>t_n} J_{31}(x)/(n\ov F(x))=0$
follows in the same way as for the first term in the last display.
The negligibility of $J_{32}(x)/(n\ov F(x))$ is again proved by
Prokhorov's inequality for iid \regvary\ \rv s,
also observing that uniformly for $x>t_n$ and fixed $\vep>0$, by Karamata's theorem  and the choice of $t_n\gg \sqrt{n\,\log n}$,
\beao
n\,|\E[\wh X]|\le n\,\E[|X|\,\1(|X|>\vep x)]\sim c\,n\,(\vep x)\P(|X|>\vep x)
=o(x)\,.
\eeao

\subsection*{The upper bound}
We start with the bound $\P(\{S_n>x\}\cap B_1)\le R_1(x)+R_2(x)$.
Since $x>t_n\gg \sqrt{n\,\log n}$ the classical Nagaev \ld\ result, \cite{nagaev:1979}, applies
to each of the $\P((S_{ij}^{(r)})'>\xi x/m)$ uniformly for $x$.
Hence, uniformly for $x>t_n$,
\beao
\dfrac{R_1(x)}{n\,\ov F(x)}\le c
\sum_{h=1}^m\P(|X_h|>\vep x\;\mid\;|X_0|>\vep x)+
c [n\P(|X|>\vep x)]=o(1)\,.
\eeao
Next, having $R_2(x) = R_{21}(x) + R_{22}(x),$ we restrict ourselves to the investigation of $R_{21}(x)$ as in the proof of Theorem~\ref{thm:main}. The relation \eqref{eq:aug22} remains true, thus we have by ${\bf RV}_3$
\beao
\frac{R_{21}(x)}{n \ov F(x)} \le c \sum_{h=1}^m\P\left(|X_h|>\vep x\;\Big|\;|X_0|>\frac{(1-\xi) x}{4m+2}\right) + c[n\ov F(x)] = o(1)
\eeao
uniformly for $x>t_n.$
\par
Next we comment on $\P(\{S_n>x\}\cap B_2)$. For any small $\delta>0$ write
$A_\delta=\cup_{i=1}^n\{X_i>(1-\delta)\,x\}$.
Thus, we derive
\beao
\P(\{S_n>x\}\cap B_2)&=&\P(\{S_n>x\}\cap B_2\cap A_\delta)+
  \P(\{S_n>x\}\cap B_2\cap A_\delta^c)\\
&\le & n\,\ov F((1-\delta)\,x)+  \sum_{i=1}^n \P\big(X_i+ \wh S_n^{(i)}>x\,,
|X_i|>\vep\,x\,,X_i\le (1-\delta)\,x\big)\\
&\le &n\,\ov F((1-\delta)\,x)+ \sum_{i=1}^n \P\big(\wh S_n^{(i)}>\delta \,x\,,
|X_i|>\vep\,x\big)=:C_1(x)+C_2(x)\,.
\eeao
Therefore,
\beao
\lim_{\nto}\sup_{x>t_n}\dfrac{C_1(x)}{n\,\ov F(x)}\le (1-\delta)^{-\alpha}\,,
\eeao
and the \rhs\ converges to 1 as $\delta\downarrow 0$.
We have by $m$-dependence
\beao
C_2(x)&\le & \sum_{i=1}^n\P\Big(\sum_{j\in [i-m,i+m], j\ne i} \wh X_j>\delta \,x/2\Big)\\
&&+\sum_{i=1}^n\P\Big(\sum_{j\le n, j\notin [i-m,i+m]} \wh X_j>\delta\,x/2\Big)\,\P(|X|>\vep x)\,.
\eeao
The sums in the first right-hand \pro ies can be bounded by $2mg(x)=2m\vep x$.
Choosing $2m\vep<\delta/2$ the first term vanishes. Writing $P_{n,i}(x)$ for the summands in the second term, we have
\beao
\dfrac{\sum_{i=1}^n P_{n,i}(x)}{n\,\ov F(x)}
&\le & c\,\vep^{-\alpha} \dfrac{1}{n}\sum_{i=1}^n\P\Big(\sum_{j\le n, j\not\in [i-m,i+m]} \wh X_j>\delta\,x/2\Big)\,.
\eeao
The \pro ies in the sum can by bounded uniformly for $x$ and $i$ by splitting the sum into $m$ sums of iid summands and then
applying Prokhorov's inequality. Moreover, this bound converges
to zero since we can choose $\vep>0$ arbitrarily small.
Thus $\sup_{x>t_n} C_2(x)/(n \ov F(x))$ is negligible as $\nto$.
\par
Finally, we bound $\P(\{S_n>x\}\cap B_3)=\P(\wh S_n>x)$. We split $\wh S_n$
into $m$ independent sums and apply the Fuk-Nagaev's inequality for $p>\alpha$.
Thus, we obtain for universal constants
$c,d>0$ only depending on $p$,
\beao
\P(\wh S_n>x) &\le& c \,n\, \E \big[|\wh X|^p\big] \,x^{-p} + \exp(- d\,x^2/ n)\,.
\eeao
We have by Karamata's theorem uniformly for $x>t_n$,
\beao
\dfrac{n\, \E \big[|\wh X/x|^p\big]}{n\,\ov F(x)} &=&
\dfrac{\E \big[|\wh X|^p\big]}{(\vep x)^p\,\P(|X|>\vep x)}\,\vep^p\,
\dfrac{\P(|X|>\vep\,x)}{\ov F(x)}\to c\, \vep^{p-\alpha}\,,\qquad \nto\,,
\eeao
and the \rhs\ vanishes as $\vep\downarrow 0$. Moreover, for large $x$
\beao
\dfrac{\exp(-d\,x^2/n)}{n\,\ov F(x)}=\exp(-d\,x^2/n -\log n +S(x))\le \exp(-0.5\, d\,x^2/n -\log n+ \alpha \log x)
\eeao
and the \rhs\ converges to zero as $\nto$
uniformly for $x>t_n\gg \sqrt{n \log n}$. \qed
}

\section*{Acknowledgements}
Thomas Mikosch's research is partly supported
by Danmarks Frie Forskningsfond Grant No 9040-00086B. Igor Rodionov's research is partly supported by the Russian Foundation for Basic Research Grant No
19-01-00090.

\end{document}